\documentclass[11pt]{amsart}

\usepackage{amsmath}
\usepackage{amssymb}
\usepackage{graphicx}

\usepackage{graphicx}%
\usepackage{multirow}%
\usepackage{amsmath,amssymb,amsfonts}%
\usepackage{amsthm}%
\usepackage{mathrsfs}%
\usepackage[title]{appendix}%
\usepackage{xcolor}%
\usepackage{textcomp}%
\usepackage{manyfoot}%
\usepackage{booktabs}%
\usepackage{algorithm}%
\usepackage{algorithmicx}%
\usepackage{algpseudocode}%
\usepackage{listings}%
\usepackage{url}

\usepackage{amsmath,amsfonts,amstext,amscd,amssymb,a4,enumerate,epsfig,psfrag}
\usepackage[latin1]{inputenc}
\usepackage{amsopn}
\usepackage{color} 

\newtheorem{theorem}{Theorem}[section]
\newtheorem{corollary}[theorem]{Corollary}

\newtheorem{proposition}[theorem]{Proposition}
\theoremstyle{definition}

\numberwithin{equation}{section}

\newcommand{\NN}{{\mathbb{N}}}

\newcommand{\cO}{{\mathcal{O}}}

\newcommand{\cB}{{\mathcal{B}}}
\newcommand{\cC}{{\mathcal{C}}}

\newcommand{\cF}{{\mathcal{F}}}

\newcommand{\cS}{{\mathcal{S}}}

\newcommand{\RR}{{\mathbb{R}}}
\newcommand{\CC}{{\mathbb{C}}}

\newcommand{\ind}{\operatorname{ind}}
\newcommand{\interior}{\operatorname{int}}

\newcommand{\Ran}{\operatorname{Ran}}
\newcommand{\Ker}{\operatorname{Ker}}

\newcommand{\sign}{{\operatorname{sign}}}

\newcommand{\bs}{{\bigskip}}

\title[Multiple solutions by inducing bifurcations]{Computing multiple solutions   from  knowledge of the critical set}

\author[O. Kaminski]{Otavio Kaminski}
\address[O. Kaminski]{Departamento de Educa\c{c}\~{a}o, IBC, Av.Pasteur 350, Rio de Janeiro, Brazil}
\email{{\tt otaviokaminski@gmail.com}}

\author[C. Tomei]{Carlos Tomei}
\address[C. Tomei]{Departmento de Matem\'{a}tica, PUC-RIO, R. Mq. S. Vicente 225, Rio de Janeiro, Brazil}
\email{\tt carlos.tomei@gmail.com}

\author[D. S. Monteiro]{Diego S. Monteiro}
\address[D. S. Monteiro]{Col\'{e}gio de Aplica\c{c}\~{a}o, UERJ, R. B. Itapagipe 96, Rio de Janeiro, Brazil}
\email{\tt diego\_smonteiro@hotmail.com}



\begin{document}

\begin{abstract}
{We explore a simple {\it geometric model} for functions between spaces of the same dimension (in infinite dimensions, we require that Jacobians  be Fredholm operators of index zero). The model combines standard results in analysis and topology associated with familiar global and local aspects. Functions are supposed to be proper on bounded sets. The model is valid for a large class of semilinear elliptic differential operators.
	
	It also provides a fruitful context for numerical analysis. 	
	For a function $F: X \to Y$ between real Banach spaces,  continuation methods to solve $F(x) = y$  may improve from  considerations about the global geometry of $F$.

	We consider three classes of examples. First we handle functions  from the Euclidean plane to itself, for which the reasoning behind the techniques is visualizable. The second, between spaces of dimension 15, is obtained by discretizing a nonlinear Sturm-Liouville problem for which special right hand sides admit abundant solutions. Finally, we  compute the six solutions of a semilinear elliptic equation $-\Delta u - f(u) = g$ studied by Solimini.}
\end{abstract}

\dedicatory{To Percy, with affection and admiration.}

\maketitle


\section{Introduction}\label{introduction}

We explore a simple {\it geometric model} for functions between spaces of the same dimension (in infinite dimensions, we require that Jacobians  be Fredholm operators of index zero). The model combines standard results in analysis and topology associated with familiar global and local aspects. Functions are supposed to be proper on bounded sets. The model is valid for a large class of semilinear differential operators and, as we shall see,  provides a fruitful context for numerical analysis.

\bigskip

The geometric model is a classification scheme --- given an appropriate (nonlinear) function $F$, it is identified to one of a collection of representatives described by the model. It is not surprising that numerics blends well with such a framework: on the one hand, it may help identifying a representative to a function, on the other, the numerics may benefit from knowledge of restrictions imposed by the classification.

\bigskip
Let $X$ and $Y$ be real Banach spaces, $U \subset X$ be an open set and   $\tilde F: U \to Y$ a continuous function. Let $\tilde \cC$ be the {\it critical set} of $\tilde F$, i.e., the set of points  $u \in U$ on which $\tilde F$ does not restrict to a local homeomorphism.
Noncritical points of $U$ are called {\it regular}. For $A \subset X$, let $\interior A$ be the interior of $A$. Under natural hypotheses, most critical points are of a special kind.  A critical point $x \in U$ is a {\it fold}\footnote{In a sense, folds are standard critical points, given their typical abundance.} of $\tilde F$ if  there are local changes of variables (i.e., homeomorphisms) centered at $x$ and $\tilde F(x)$ converting $\tilde F$ into
\[ (t, z) \in \RR_0 \times Z_0 \subset \RR \times Z \mapsto (t^2, z) \in \RR \times Z\]
for a real Banach space $Z$ and appropriate open neighborhoods $\RR_0$ and $Z_0$ of the origin of $\RR$ and $Z$. A point $y \in F(\cC) \subset Y$ is a {\it standard critical value} if  its preimages contain regular points and exactly one fold point.
For a closed set $D \subset U$, consider the restriction $F: D \to Y$ and define its critical set $\cC$ to be $\tilde \cC \cap D$. In particular, $\partial D$ may contain critical and regular points of $F$.

A point $x \in \partial D$ is a {\it standard boundary point} if there is an open neighborhood  $U_x$ of $U \subset X$ on which $F$ restricts to a local homeomorphism $\tilde F$ which, after changes of variables, becomes  
\[ \tilde \Phi: (t, z) \in \RR_0 \times Z_0 \mapsto (t,z) \in \RR_0 \times Z_0\] 
where $\RR_0$ and $Z_0$ are as above, and
\[\tilde \Psi(\partial D \cap  U_x) = (0, Z_0) \ \hbox{and} \ \tilde \Psi(D \cap U) = \{ (t, Z_0), t \ge 0\}. \]
A point  $y \in F(\partial D) \subset Y$ is a {\it standard boundary value} if its preimages contain only regular points in $\interior D$ together with an additional single point  $x \in \partial D$ which is a standard boundary point.

Let $A \subset X, B \subset Y$. A continuous function $F: A \to B$ is {\it proper} if the inverse of a compact set in $B$ is a compact set in $A$. Continuous functions between finite dimensional spaces are automatically compact if defined on  bounded sets. 
Properness in infinite dimensions is more elusive, as there are closed, bounded sets which are not compact. Still, there is an interesting class of  elliptic differential operators for which this is not an issue: more generally, there are functions $F: A \to B$ for which the inverse of a compact set of $B$, if bounded, is compact in $A$.  We provide simple, familiar examples in Section \ref{primeiroexemplo}.

For a continuous function $F: D \subset X \to Y$ with critical set $\cC$, define the {\it flower} $\mathcal{F} = F^{-1}(F(\cC \cup \partial D)) \subset D$. 
Domain  $D$ and codomain $Y$   split in {\it tiles}, the connected components of the complements $\cF^c \subset D$ and $F^c(\cC \cup \partial D) \subset Y$   respectively. By connectivity, a tile in the codomain either does not intersect or it is contained in the image $F(D)$ and we refer to it as an {\it image tile}.

\medskip

\begin{theorem} [Geometric model] \label{theo:model} Let $X, Y, U, D$ as above and $F: D \to Y$ be a proper function.
	
	\begin{enumerate}
		\item The sets $\cC, \partial D, F(\cC \cup \partial D)$ and $\cF = F^{-1} (F(\cC \cup \partial D))$ are closed. Domain (resp. image) tiles are open sets in $D$ (resp. $Y$).
		\item 	The restriction of $F$ to a domain tile is a covering map of finite degree onto an image tile. In particular, all points in an image tile have the same (finite) number of preimages.
		\item $F$ takes domain tiles sharing a fold in their boundary  to the same image tile. 
		\item 
		The number of preimages of two image tiles sharing a standard boundary value (resp. a standard critical value) in their boundary differs by 1 (resp. by 2). 
	\end{enumerate}
\end{theorem}

\medskip

The text handles a number of examples of such functions, ranging from functions from the plane to the plane to nonlinear differential operators. We alternate numerical and theoretical material, aiming at a self-contained presentation.

\bigskip
We now relate these geometric concepts to the general problem of computing a solution $u_1 \in D$ of $F(u_1)=g_1$. Most continuation methods belong to two broad categories (excellent references are \cite{ALLGOWER,RHEINBOLDT, UECKER}). In the first, given some  $u_0 \in D$ for which $F(u_0) = g_0$, one  proceeds to invert a curve $g(t)$ joining $g(0) = g_0$ and $g(1) = g_1$ (usually a segment), starting at $t=0$. The implementation of the inversion is not relevant at this point, but its realization (i.e., completing inversion up to $t=1$) is of essence.

In the second category, for which frequently $g_1 = 0$,  one first extends $F$ to a function $G(u,t): D \times \RR \to Y$ for which a simple curve $d = (u(t),t)$ in the domain  satisfies $G(u(0),0) = g_1$. Assuming differentiability, one searches for additional preimages by first identifying {\it bifurcation points} $u(t_c)$, for which $DG(u(t_c), t_c): X \times \RR \to Y$ is not surjective for $t_c \in [0,1]$.
Such points generically yield new branches of solutions $(\tilde u(t),t)$: if they extend to $t=1$, new preimages of $g_1$ are obtained. Additional bifurcations may arise along such branches. Usually, the specification of the curve $d = (u(t),t)$ is restricted, being related to the symbolic expression of $F$.

Clearly, the critical set $\cC$ of $F$ and its image $F(\cC)$ are fundamental objects in such continuation methods. In the first category, from the geometric model, if the segment $g(t)$ lies within an image tile, a connected component of $Y \setminus F(\cC)$, inversion is feasible at least on theoretical grounds, by the unique lifting property of covering maps \cite{HATCHER}. In the second, bifurcation points of $G$ belong to  the set of critical points of the part of $DG$ related to the derivative in $u$. It is rather surprising that there is no reference to $\cC$ in the choice of the curves that guide both procedures. 

More than an algorithm, we suggest a simple strategy. A data bank of points with abundant, computable preimages would be helpful in the solution of the original task, solving $F(u) = g$. Indeed, points $g_0$ for which many preimages are known give more opportunities for complete inversion using methods in the first category.  We search for such points by {\it inducing bifurcations}. From scarse knowledge of the critical set $\cC$ of  $F$, we choose curves  $c \subset D$ through a point $u_0$ with substantial intersection with the critical set $\cC$, and compute the associated {\it bifurcation diagram}, the connected component of $F^{-1}(F(c))$ containing $u_0$. 
The strategy requires that we may decide  if $u \in \cC$ or  $u \in\partial D$. Inversion near critical points relies on spectral information, in the spirit of Section 3 of \cite{UECKER} (see also Appendix \ref{continuationatafold}).

\bigskip

How does the strategy compare with more familiar algorithms?
The freedom in the choice of the curve $c$ in the domain allows to  generically perform continuation only at regular points and folds. Moreover, by choosing both $u_0$ and $c$ appropriately, we are freer to search for situations with abundant bifurcations. 
For a more pragmatic evaluation of the strategy, we present it in action in three different contexts, which we describe superficially below.

Continuation methods received a recent boost from ideas by Farrell, Birkisson and Funke \cite{FARRELL}, which improved  an elegant deflation strategy originally suggested by Brown and Gearhart \cite{BROWN}. The techniques can be grafted to our scheme, but we do not consider the issue.

Our approach started with the study of  proper functions $F: \RR^2 \to \RR^2 $  by  Malta, Saldanha and Tomei in the late eighties (\cite{MST1}), leading to $2\times 2$, a software that computes preimages of $F$, together with other relevant geometrical objects\footnote{Available at \url{www.im-uff.mat.br/puc-rio/2x2}. Special thanks to Humberto Bortolossi! }. Under generic conditions, the authors obtain a characterization of critical sets: given finite sets of curves $\{\cC_i\}$ and $\{\cS_i\}$, and finite points $\{p_{ik} \in \cC_i\}$, one can decide if there is a proper, generic function $F$ whose critical set consists of the curves $C_i$, with images $ \cS_i= F(\cC_i) $ and cusps\footnote{Informally, for functions in the plane, cusps are the second most frequent critical points \cite{Whitney}.}\footnote{Notice that the theorem does not claim the existence or not of additional critical curves. Indeed, this cannot be verified in a computational environment on which only a finite number of function evaluations are permitted, as opposed to, for example,  handling estimates on derivatives. In particular, in the usual environment, given $F$,  no computer program may decide if all the solutions of $F(u) = g$ are known.} at the points $\{p_{ik}\}$. Given a function $F$, $2 \times 2$ obtains some critical curves $\mathcal{C}_i$, its images $\cS_i$ and their cusps $\{p_{ik}\}$: if the characterization does not hold, it provides the program with information about where to search for additional critical curves.

Two important features of the 2-dimensional context do not extend:
(a) a description of the critical set as a list of critical points,
(b) the identification of higher order singularities.
A realistic implementation for $n > 2$ led us to the current text. Our research was  strongly motivated by the underlying geometry of the celebrated Ambrosetti-Prodi theorem (\cite{AP,MM}). Subsequent articles (\cite{AP,BERGER,CALNETO,SMILEY,KAMINSKI}) provided  information about the geometry of semilinear elliptic operators, with  implications to the underlying numerics. 

The geometric approach in $2 \times 2$ found application in different scenarios (for ODEs,
\cite{TOMEIEBUENO,BURGHELEA,TELES}; for PDEs, \cite{CALNETO,KAMINSKI}).

\bigskip
In Section \ref{model} we prove Theorem \ref{theo:model} and in Section \ref{Visual}   we apply  the strategy to some visualizable examples. There are strong conceptual affinities among the examples along the text, as we shall see.

The second class of examples, discussed in Section \ref{SturmLiouville},  arises from the discretization $F^h$ of a nonlinear Sturm-Liouville problem,
$$ F(u) =  - u'' + f(u) = g  , \quad u(0) = 0 = u(\pi)  .$$
The discretized equation $F^h(u)= g$ for $u , g \in \RR^n$ has an unexpectedly high number of solutions (\cite{TELES}), and most do not admit a continuous limit, but we consider the discretized problem as an example of interest by itself. 
In \cite{ALLGOWER2}, Allgower, Cruceanu and Tavener approached numerical solvability of semilinear elliptic equations by first obtaining good approximations of solutions from discretized versions of the problem. A filtering strategy eliminates some  discrete solutions which do not yield a continuous limit. As they remark, results on the number of solutions are abundant, but not really precise. 

Geometry intervenes as follows (\cite{TELES,TOMEIEBUENO, BURGHELEA, BURGHELEA2}). The connected components $\cC_i, i=1, \ldots, n$, of the critical set of $F^h$ are graphs of functions $\gamma_i: V^\perp \to V$, where $V$ is the one dimensional space generated by the ground state of the standard discretization of the  $u \mapsto -u''$ with Dirichlet boundary conditions, which, as is well known, is the evaluation of  $u(x) = \sin x$ at points of a uniform mesh on $[0, \pi]$.  Curves in the domain with abundant critical points are easy to identify:  follow straight lines of the form $u_0 + t\sin (x)$, for different values of $u_0$. Parallel lines frequently give rise to disjoint sets of solutions.

The third context, equation (\ref{AP}) in Section \ref{Solimini}, is a semilinear operator $F(u) = -\Delta u - f(u)$ considered by Solimini (\cite{SOLIMINI}) acting on functions satisfying Dirichlet boundary conditions in an annulus. A special right hand side $g$ has exactly six preimages, which are obtained by inverting  lines of the form $u_0 + t \phi_0$, where $\phi_0$ is the ground state of the free (Dirichlet) Laplacian. Minor perturbations of such lines circumvent nongeneric difficulties. 

The visual examples in the plane  in a sense are harder than the discrete and the continuous functions arising from semilinear elliptic equations. Indeed, in the last two cases all we have to know from the critical set is the fact that it is intercepted abundantly by lines which are parallel to the ground state of the free operator, a simple consequence of a min-max argument.

Section \ref{terminology} is an Appendix in which we collect some basic definitions and properties for the reader's convenience.
In Section \ref{continuationatafold} we handle inversion of segments in the neighborhood of the critical set $\cC$. In particular, we avoid the difficulty of partitioning Jacobians for discretizations of infinite dimensional problems. Instead of using variables  related to arc length \cite{ALLGOWER},
we work with spectral variables, which are especially convenient for the strategy.

We used $2 \times 2$ and MATLAB for graphs and numerical routines.

\medskip
\noindent{\bf Acknowledgements:} Tomei was partly supported by grants from Stonelab, FAPERJ E-26/200-980/2022 and CNPq 304742/2021-0. Kaminski and Moreira acknowledge the support from CAPES and CNPq. The authors are very grateful to a referee, which provided a careful, generous reading.

\section{Basic geometry} \label{model}

With the notation of the Introduction, we consider functions $F: D \subset U \subset X \to Y$ between Banach spaces $X$ and $Y$, for which we defined critical and regular points and domain and image tiles.

{\it In the figures, a point $A$ in the image of $F$ has preimages $A_i$: $A = F(A_i)$.}


\subsection{A simple function satisfying the geometric model}

Consider the proper function
$$ F:D = U = X = \RR^2 \to Y =  \RR^2 \ , \quad (x,y) \mapsto  (x^2 - y^2 +x, 2 xy - y)\ . $$

\begin{figure} [ht] 
	\begin{centering}
		\includegraphics[scale=0.33]{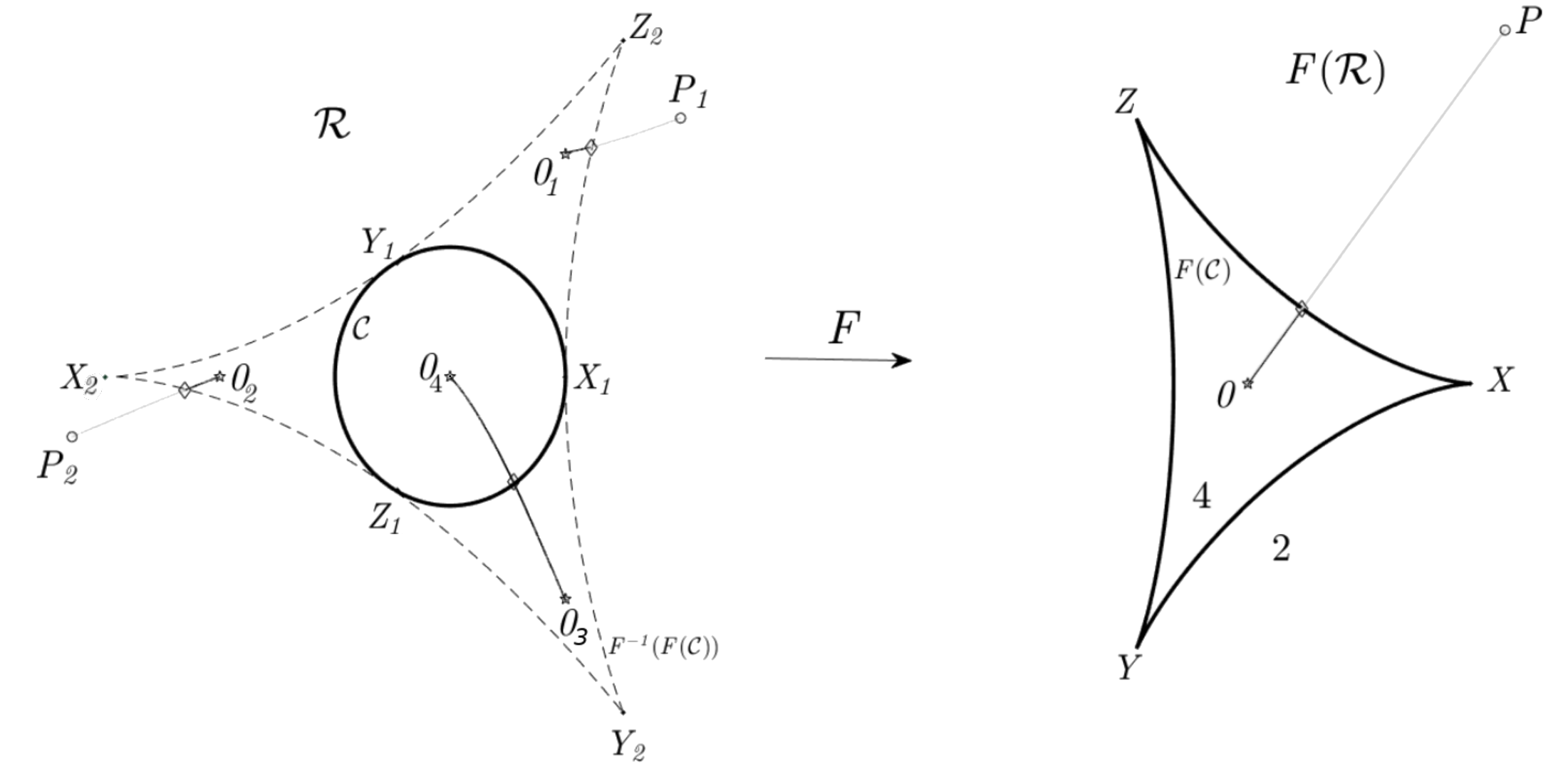}
		\caption{Five domain tiles, two image tiles. }
		\label{fig:z2zbarra}
	\end{centering}
\end{figure}

Its domain
on the left of Figure \ref{fig:z2zbarra}, contains the critical set $\cC$, a circle, and the flower $\mathcal{F}$, consisting of  triangles\footnote{By a triangle, we mean the region bounded by three arcs; the boundary of the disk is a triangle too.} with vertices $X_i, Y_i$ and $Z_i, i=1,2$. There are five domain tiles and two  image tiles, one bounded and one unbounded. The function $F$ is a homeomorphism from each of the four bounded domain tiles  to the bounded image tile. Thus, points in the bounded image tile have four preimages. In particular, $F$ has roots $0_1, \ldots, 0_4$, the preimages of $0$. The unbounded tile $\mathcal{R}$ in the domain is taken to  the unbounded tile $F(\mathcal{R})$, but the function is not bijective: each point in  $F(\mathcal{R})$ has {\it two} preimages, both in $\mathcal{R}$. More geometrically, all restrictions of $F$ to tiles are covering maps (\cite{HATCHER}), and thus must be diffeomorphisms when their image is simply connected. In Figure \ref{fig:z2zbarra}, the numbers on the image tiles are the (constant) number of preimages of points in each tile.

\subsection{The simplest critical points: folds} \label{basicfolds}

In Figure \ref{fig:z2zbarra}, points in different tiles in the image (necessarily adjacent, in this case, i.e., sharing an arc of images of critical points) have their number of preimages differing by two.
Inversion of the segment $[(0,0) , P]$ by continuation gives rise to two subsegments, whose interiors have two and four preimages. Starting from $P$ with initial condition $P_1$ or $P_2$, inversion carries through to $0$ without difficulties, giving rise to roots $0_1$ and $0_2$. However, when inverting from $0$ with initial conditions $0_3$ and $0_4$, continuation is interrupted. This is the expected from the behavior of the function $F$ at a fold, as we now outline.

\begin{figure} [ht] 
	\begin{centering}
		\includegraphics[scale=0.3]{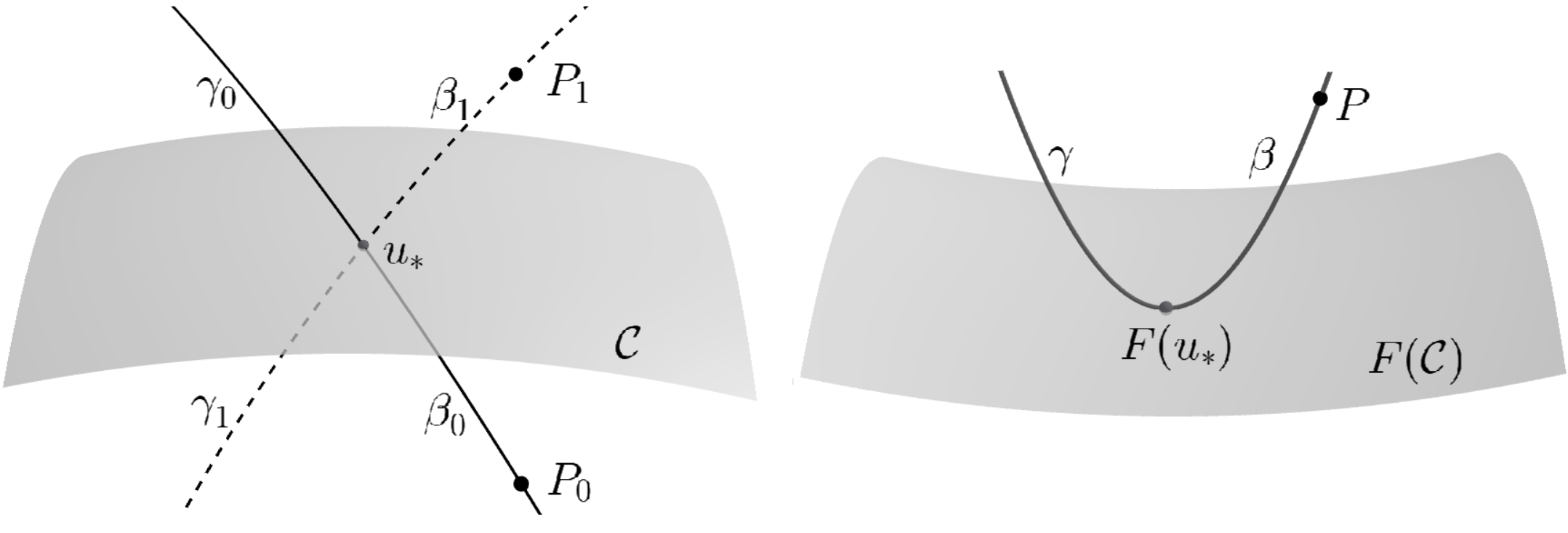} 
		\caption{Near a fold $u_\ast$.}
		\label{fig:extensao2}
	\end{centering}
\end{figure}

A critical point $x \in U$ is a {\it fold} of $F$ if and only if there are local changes of variables centered at $x$ and $F(x)$ converting $F$ into
\[ (t, z) \in \RR \times Z \mapsto (t^2, z) \in \RR \times Z\]
for some real Banach space $Z$.  On Figure \ref{fig:extensao2}, point $u_\ast$ is a fold. The vertical segment through $u_\ast$ splits in two segments, $\beta_0$ and $\beta_2$. The inverse of $F(\beta_0)$ yields again $\beta_0$ and another arc $\beta_1$, also shown in Figure \ref{fig:extensao2}. In a similar fashion, the inverse of $F(\beta_2)$ contains $\beta_2$ and another arc $\beta_3$. Notice the similarity with a bifurcation diagram: two new branches $\beta_1$ and $\beta_3$ emanate from $u_*$.

\medskip
For a characterization of folds in Banach spaces, see Appendix \ref{folds}.

\subsection{Boundary points}

\begin{figure} [ht] 
	\begin{centering}
				\includegraphics[scale=0.35]{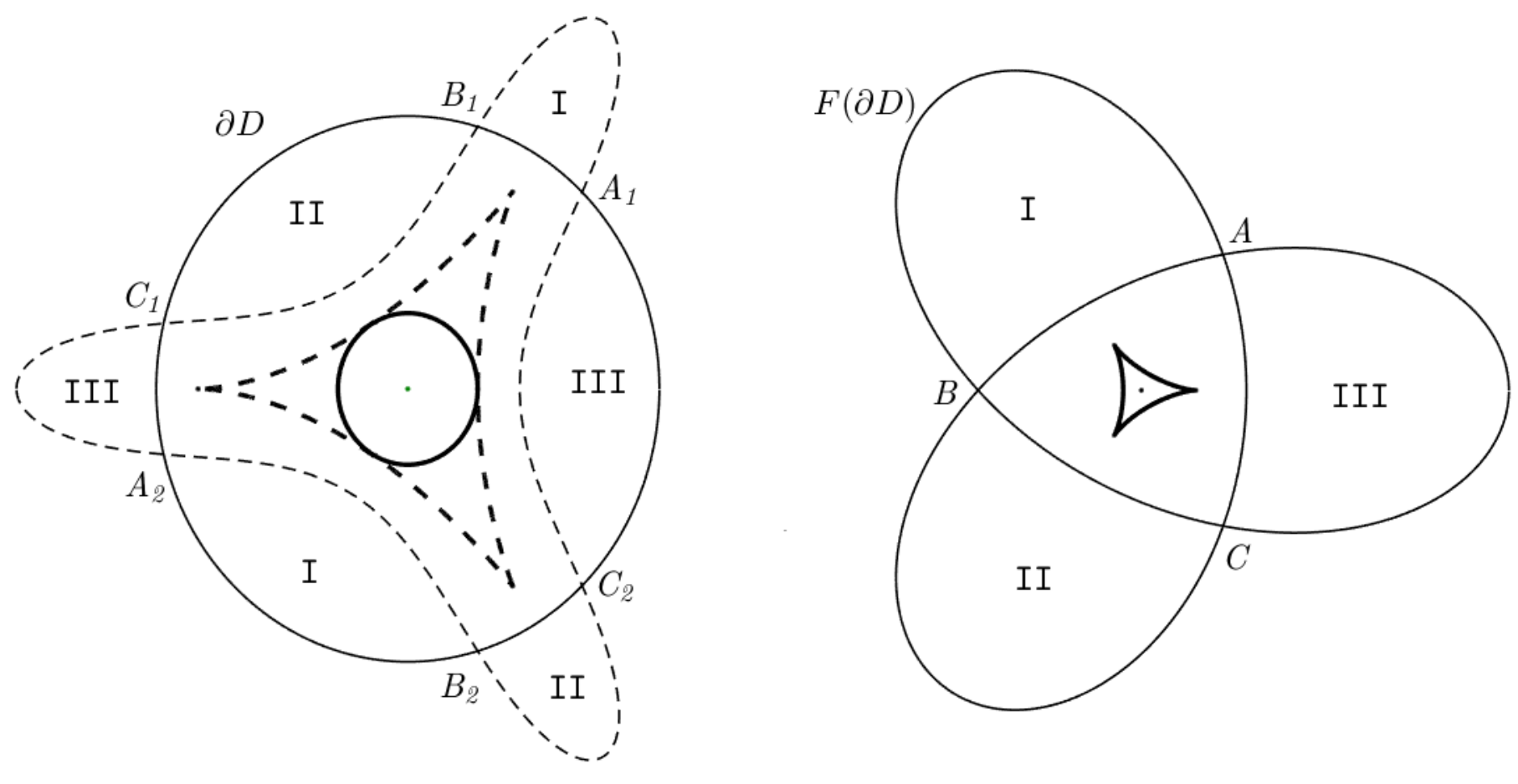}
		\caption{Behavior at a well behaved boundary. }
		\label{fig:fronteira}
	\end{centering}
\end{figure}

We consider an expressive example. Restrict $F$ above to a disk $D$ with boundary $\partial D$ described in Figure \ref{fig:fronteira}. The flower now also includes $\partial D$ and an additional dotted curve $F^{-1} (F(\partial D))$, containing the preimages $A_i, B_i, C_i, i=1,2$ of points $A, B$ and $C$. Each one of the image tiles I, II and III has a single preimage  inside $D$, the other being outside. All points in $F(\partial D)$ with the exception of $A$, $B$ and $C$ are standard boundary points Trespassing an arc of points of $F(\partial D)$ at a standard boundary point only changes the number of preimages by one.

\subsection{Proper functions} \label{primeiroexemplo}

A continuous function $F: D \to Y$ is {\it proper} if the inverse of a compact set of $Y$ is a compact set of $D$. More generally, $F$ is {\it proper on bounded sets} (b-proper) if its restriction $F_B: B \subset D \to Y$ to bounded, closed sets $B$ is proper.

A continuous function $F: \RR^n \to \RR^n$ is proper if and only if $ \| F(x) \| = \infty$ as $\|x\| \to \infty$. Many elliptic semilinear operators are b-proper: we provide a standard example in Corollary \ref{elipticos}.

\medskip

\begin{proposition} \label{gerais} Let $Y$ be a real Banach space and $G: Y \to Y$, $G(u) = u + \Phi(u)$, be continuous. Suppose that, for any closed ball $B \subset Y$, $\overline{\Phi(B)}$ is a compact set. Then $G$ is b-proper. Moreover,
	$G$ is proper if and only if $\| G(x) \| \to \infty$ as $\|x\| \to \infty$.
\end{proposition}

\begin{proof} We prove that $G$ is b-proper. 
	Let $B$ be a closed, bounded subset of $Y$: we show that the restriction $\hat G: B \to Y$ of $G$ is proper. Let $K \subset Y$ be a compact set. Take
	a sequence $y_n \in K$ with $y_n\to y_\infty \in K$, and $u_n \in B \subset Y$ such that $G(u_n) = u_n + \Phi(u_n) = y_n$. For an appropriate subsequence, $\Phi(u_{n_m}) \to w \in Y$ and then $u_{n_m}= - \Phi(u_{n_m}) + y_{n_m}\to -w + z_\infty \in Y$. Thus $B \cap G^{-1}(K)$ is compact and $G$ is b-proper. The remaining claim is left to the reader.
\end{proof}

\begin{corollary} \label{elipticos} Let $\Omega \subset \RR^n$ be a bounded set with smooth boundary. For a smooth $f: \RR\to \RR$, set $F:X = C^{2,\alpha}_D(\Omega) \to Y = C^{0,\alpha}(\Omega)$ given by $F(u) = - \Delta u + f(u)$. Then $F$ is b-proper.
\end{corollary}

\medskip
Here, $C^{2,\alpha}_D(\Omega) $ is the  H\"older space of functions equal to zero on $\partial \Omega$, $\alpha \in (0,1)$.



\begin{proof} Set $ G: Y \to Y, \ G(v) = v + \Phi(v)$, where $\Phi(v) = f(-\Delta^{-1} v)$. Recall that $-\Delta: X \to Y$ is an isomorphism and the inclusion $\iota: X \to Y $ is compact. Thus, for a bounded ball $B \subset Y$, the set $\overline{\iota(-\Delta^{-1} B)}$ is compact. From Theorem 2.1 in \cite{CHIAPPINELLI}, the Nemitski map $F: Y \to Y, u \mapsto f(u) $ is continuous on bounded sets and thus $\Phi(v)$ satisfies the hypotheses of the proposition above. 
\end{proof}

\bs

\subsection{Proof of Theorem \ref{theo:model}}

Denote the complement of a set $A \subset B$  by $A^c \subset B$.
We start with item (1).
Clearly $\cC$ and $\partial D$ are closed sets, being complements of open sets. It is easy to see that the image of closed set under a proper map is also closed, thus $F(\cC)$ is closed, and by the continuity of $F$, the flower $\cF$ also is. By taking complements, tiles are open sets, proving (1).

By definition, a domain tile $T_D$ is an open, connected set. As $F$ is proper, $F(T_D)$ also is and must lie in some image tile $T_C$. 
We show by contradiction that $T_C = F(T_D)$. Indeed, suppose $y \in T_C \setminus F(T_D)$ with $y_n \in F(T_D)$ such that $y_n\to y$.	As in the previous paragraph, from the properness of $F: D \to Y$ and the compactness of $\{y_n, n \in \NN\} \cup \{y\}$, there is a convergent subsequence $x_{n_m} \in T_D, x_{n_m} \to x \in \overline{T_D}$ for which $F(x_{n_m}) = y_{n_m} \to y$ and thus $F(x) = y$.
There are two cases to consider. If $x \in T_D$, then $y \in F(T_D)$, a contradiction.
Otherwise, $x \in \partial T_D \subset \cF = F^{-1} (F(\cC \cup \partial D))$,  contradicting $F(x) = y \in T_C \subset F^c(\cC \cup \partial D)$. Thus the restriction $\tilde F: T_D \to T_C$ is surjective and a local homeomorphism at each point of its domain, as $T_D$ consists of regular points. From the argument, we also have that $F( \partial T_D) \subset \partial T_C$, which implies that the restriction  $\tilde F: T_D \to T_C$ is proper. 

Each point of $T_C$ has finitely many preimages in $T_D$. Indeed, as in the arguments above, infinitely many preimages of $y \in T_C$ would accumulate at some point  $x \in T_D$ by the properness of $\tilde F$, contradicting the fact that $\tilde F$ is a local homeomorphism at $x$. In order to prove that $\tilde F$ is a covering map, one must find an open neighborhood $V_y$ of each point $y \in T_C$ for which there are disjoint neighborhoods of $F^{-1}(V_y)$ centered at the preimages $x_i$ of $y$. This is trivial, as the set $\{ x_i\}$ is finite and $\tilde F$ is a local homeomorphism. 

Every covering map satisfies the second claim in (2): we give details for the reader's convenience. we give details for the reader's convenience. By connectivity, it suffices to show that points sufficiently close to $y \in T_C$ have the same number of preimages.
Suppose then $F^{-1}(y) = \{x_1, \ldots, x_k\} \subset T_D$, a collection of regular points: there must be sufficiently small, non-intersecting neighborhoods $V_{x_i}, i=1, \ldots,k$ of the points $x_i$ and $V_y \subset T$, such that the restrictions $F:V_{x_i} \to V_y$ are homeomorphisms. Thus, points in $V_y$ have at least $k$ preimages. Assume by contradiction  a sequence $y_n \to y, y_n \in V_y$, such that each  $y_n$ admits  an additional preimage $x_n^\ast$, necessarily outside of $\cup_i V_{x_i}$. By properness, for a convergent subsequence $x_{m_n}^\ast \to x_\infty^\ast$, $F(x_\infty^\ast)= y$, and $x_\infty^\ast \ne x_i, i=1,\ldots,k$. This completes the proof of (2).

Item (3) is a consequence of the local form of a fold $x \in \partial T_{D_1} \cap \partial T_{D_2}$. For nearby points $x_1 \in T_{D_1}$ and $T_{D_2}$ for which $F(x_1) = F(x_2) \in F_C$ we must have $F(D_1) = F(D_2)$ from the arguments above.

Finally item (4) combines the local form of $F$ at folds and standard boundary points with the proofs of the items above. Theorem \ref{theo:model} is now proved.

\bigskip
The result admits natural extensions. The image of the function $F: \RR^2 \to \RR^2$ given by $F(x,y) = (x^2, y^2)$ is the positive quadrant. Leaving it through a boundary value  different from the origin implies a change of number of preimages equal to four: this is because such boundary value has two preimages which are folds.

\bigskip
We avoid using the adjective generic, which is frequently endowed with a specific meaning by some form of Sard's theorem, or more generally, some transversality conditions coupled with Baire's category theorem. In most examples of functions $F: X \to Y$, one naturally satisfies additional hypotheses. At a risk of sounding pedantic, 
examples are the following.

\smallskip \noindent (H1) $\cC, \partial D, F(\cC)$ and $F(\partial D)$ have empty interior.

\smallskip \noindent (H2) A dense subset of $\cC$ consists of folds.

\smallskip

Still,  in Section \ref{table} we present a (finite dimensional) situation for which $F$ may have a critical set with nonempty interior. In this text, take the standard approach in numerical analysis: one proceeds (carefully, consciously) with the inversion process and accepts an occasional breakdown.

\section{The strategy on some visualizable applications} \label{Visual}

Following the strategy presented in the Introduction, given a function $F: D \subset X \to Y$, we search for points $y \in Y$ with a large number of preimages $x_i$. 
From scarse information about the critical set $\cC$ of $F$, we choose lines $c$ on the domain passing through a point $x_0$ which intersect $\cC$ abundantly. One should then expect that the {\it bifurcation diagram} $\cB$, the connected component of $F^{-1}(F(c)) \cap D$ containing $x_0$, contains some preimages of $x_0$. Implementation of inversion  follows the standard continuation strategies (\cite{ALLGOWER,RHEINBOLDT, UECKER}; see also Appendix \ref{continuationatafold}).

\begin{figure} [ht] 
	\begin{centering}
		\includegraphics[scale=0.30]{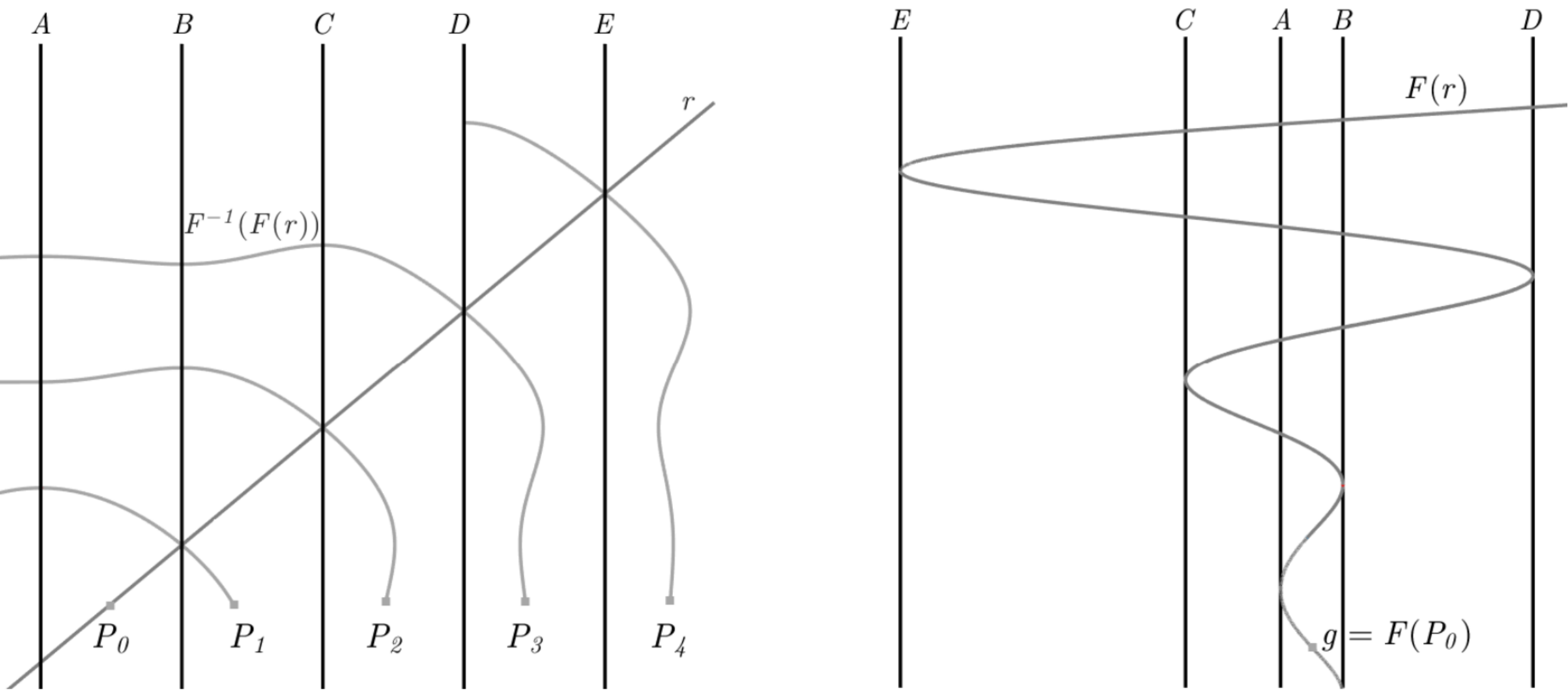}
		\caption{Pleats along fold lines $A, B, C,...$, the bifurcation diagram $\mathcal{B}$ and  some preimages $P_i$ of $g$.  }
		\label{figbabado}
	\end{centering}
\end{figure}

The image of the function
$$ F: (-0.1, \infty) \times \RR \to  \mathbb{R}^2 \ , \quad (x,y) \mapsto (\cos (x) - x^2\cos (x) + 2x \sin (x), y)$$
is well represented by  a piece of cloth pleated  along vertical lines.
As indicated in Figure \ref{figbabado} for $x$ restricted to an interval, the critical set $\mathcal{C}=\{(k\pi,y), \  k \in \NN , \ y \in \mathbb{R}\}$ start with vertical lines $A, B, C...$ and its image $F(\cC)$ also consists of vertical lines: critical points are folds.   

From the explicit knowledge of $\cC$,  a curve $c$ with abundant critical points is easy to identify: take  a straight line $r$ as in Figure \ref{figbabado}. The image $F(r)$ oscillates among images of critical lines. Given any regular point $P_0$, inversion of $F(r)$ yields $\mathcal{B}$, the bifurcation diagram $\mathcal{B}$ associated with $r$ from $P_0$, the connected component of $F^{-1}(F(r))$ containing $P_0$.  Bifurcation points, at the intersection of $r$ and $\mathcal{C}$, give rise to other preimages $P_i$ of $g=F(P_0)$. There are infinitely many such preimages.

\begin{figure} [ht] 
	\begin{centering}
		\includegraphics[height=150pt,width=370pt]{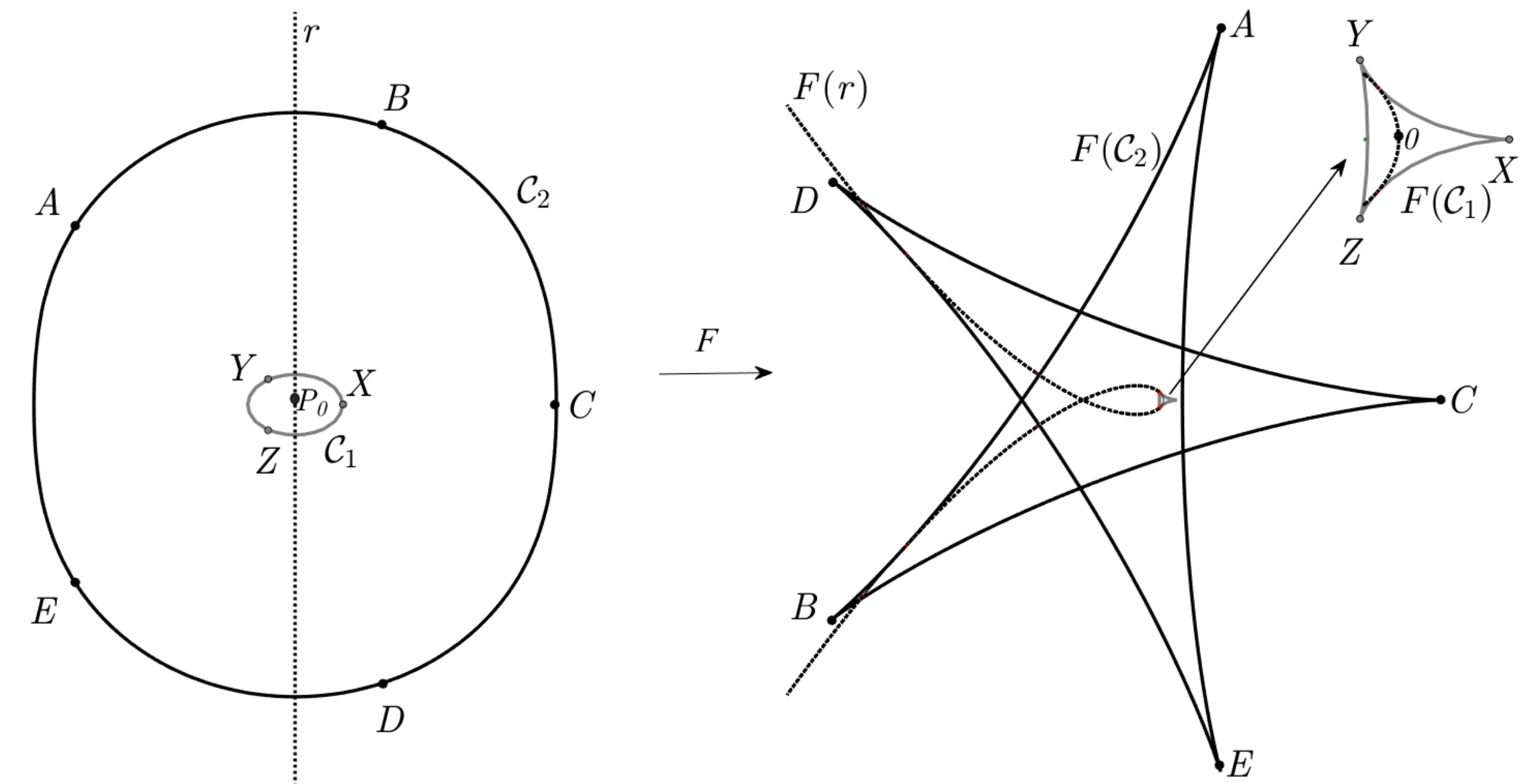}
		\caption{The critical set $\cC$, the line $r$ and their images.}
		\label{fig:estrela}
	\end{centering}
\end{figure}

\bigskip
Consider now the smooth (not analytic) function
$$ F:\CC \to \CC \ , \quad z \mapsto  z^3 +\frac{5}{2}\ {\overline{z}}^2+z \ . $$
The critical set $\cC$ consists of the two curves $\mathcal{C} _1$  and $\mathcal{C} _2$, shown in Figure \ref{fig:estrela}, which roughly bound the three different regimes of the function, $z \sim 0, 1$ and $\infty$, where $F$ behaves  like $z$, $ {\overline{z}}^2$ and $z^3$ respectively. The  three  regimes already suggest that a line through the origin must hit the critical set at least four times.

\begin{figure} [ht] 
	\begin{centering}
				\includegraphics[height=150pt,width=370pt]{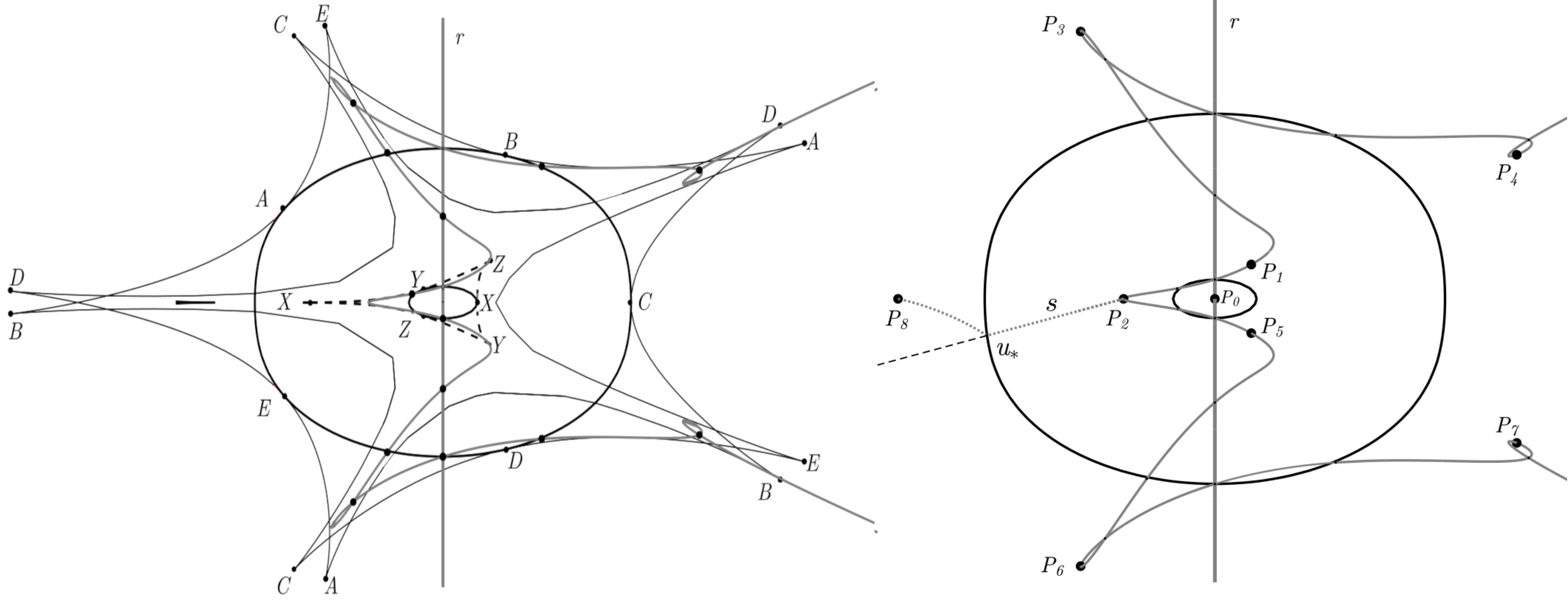}
		\caption{On the left, $\cC, r,  \mathcal{F}$ and $\mathcal{B}$. On the right, $r, \cC, \mathcal{B}$ and an extension yielding a zero $P_8$.}
		\label{fig:novinha}
	\end{centering}
\end{figure}

We count preimages with Theorem  \ref{theo:model}.  From its behavior at infinity, $F$ is proper. Points in the unbounded image tile have three preimages, as for $z \sim \infty$, the function is cubic. The  unbounded tile in the domain covers the unbounded tile in the image of the right hand side  three times. Each of the five spiked image tiles has five preimages, and the annulus surrounding the small triangle $XYZ$, seven. Finally, the interior of the small triangle $XYZ$ has nine preimages: in particular, $F$ has  nine zeros. The flower,  in Figure \ref{fig:novinha}, illustrates these facts.

Let $r$ be the vertical axis, $P_0=(0,0) \in r$. Figure \ref{fig:estrela} shows $r$ and $F(r)$ and
Figure \ref{fig:novinha}  the  flower $ \mathcal{F} = F^{-1}(F(\mathcal{C}_1)) \cup F^{-1}(F(\mathcal{C}_2))$: dotted black lines and $\cC_1$ form  $F^{-1}(F(\mathcal{C}_1))$, while $F^{-1}(F(\mathcal{C}_2))$ consists of continuous black lines.  Amplification of $ \mathcal{F}$ shows  five thin triangles, one in each `petal' (one is visible in the petal on the left), each a full preimage of the triangle $XYZ$ in the image (the enlarged detail in Figure 	\ref{fig:estrela}). Add the other four preimages at tiles bounded by $F^{-1}(F(\mathcal{C}_1))$ to spot the nine zeros of $F$.  The flower is computationally expensive and is not necessary for the strategy we present.

On the left of Figure \ref{fig:novinha} are shown $\cC$, $\mathcal{F}$, $r$ and $\mathcal{B}$, the bifurcation diagram associated with $r$ from $P_0=(0,0)$. To emphasize $r$ and $\mathcal{B}$, we removed $\mathcal{F}$ on the right of Figure \ref{fig:novinha}. The sets $r$ and  $\mathcal{C}$ meet at four folds. The set $\mathcal{B}$ contains the line $r$ and eight of the nine zeros of $F$. The missing zero, $P_8$, is on the petal on the left of Figure \ref{fig:novinha}.


The  branches originating from the four critical points in $r$ lead by continuation to additional zeros $P_i$ of $F$.  What about the missing zero? For the half-line $s$ joining $P_2$ to infinity in  Figure \ref{fig:novinha}, bifurcation at $u_*\in F^{-1}(F(s))$  obtains $P_8$. For completeness,
$$
\begin{array}{c}
P_1=( 0.2141, 0.3313) \ , \quad P_2 = (-0.5367,0.0000)\ , \quad P_3=( -0.7893, 2.5802)\ , \\
\quad P_4=(1.7752,1.3903)\ ,  \quad P_5=( 0.2141,-0.3313)\ , \quad P_6=(-0.7893, -2.5802)\ ,\\
P_7=(1.7752,-1.3903)\ , \quad P_8=(-1.8633,0.0000) \ .
\end{array}
$$

Why choose the half-line $s$?  Informally, $P_2$ is the only zero of $F$ whose neighboring zeros in $\mathcal{B}$ are farther from infinity, in the sense that it takes more critical crossings to get from them to the unbounded component. There is some luck in the process, but any line starting from the origin roughly along the direction of $s$ would also yield $P_8$. We return to this issue in the next section.

\subsection{Singularities and global properties of $F$} \label{winding}

\medskip
Curves $r$ with different intersection with the same critical component of $\cC$ may yield different zeros. Indeed, this is the case for the three different half-lines from 0 in the second example of the previous section. This happens also for the simpler function $F$ in Figure \ref{fig:z2zbarra}.  Informally, $F$ at a fold point $p$ of $\cC$ looks like a mirror, in the sense that points on both sides of the arc $\cC_i$ of $\cC$ near $p$ are taken to the same side of the arc $F(\cC_i)$ near $F(p)$, as shown in Figure \ref{fig:extensao2}. In the example in Figure \ref{fig:z2zbarra}, at the three points  which are not folds, one might think of adjacent broken mirrors.  Close to the image $X = F(X_1)$ of such points, say $X_1$, there are points with three preimages close to $X_1$. 

In higher dimensions (in particular, infinite dimensional spaces), there are critical points at which arbitrarily many broken mirrors coalesce. At such higher order singularities, there are points with clusters of preimages.  H. McKean  conjectured the existence  of arbitrarily deep singularities for the operator $F(u) = - u'' + u^2$, where $u$ satisfies Dirichlet conditions in $[0,1]$, and the result was proved in \cite{ARDILA}.

The original algorithm, $2 \times 2$, identifies and explore cusps in the two dimensional context. In more dimensions, the computation of higher order singularities is essentially unfeasible: instead, the strategy in this text must perform more line searches, by choosing curves $r$ intercepting the critical set at different mirrors.

Abundant preimages are also related to the fact that  there may be tiles on which the covering map induced by the restriction $F$  is not injective, as in the case of the annulus tile in the domain of $F$ of this section. This in turn requires the image tile to have nontrivial topology (more precisely, a nontrivial fundamental group), as simply connected domains are covered only by homeomorphisms.

\section{Discretized nonlinear Sturm-Liouville operators} \label{SturmLiouville}

We consider variations of the nonlinear Sturm-Liouville equation
\begin{equation} \label{SturmLiouville}
F(u) =  -u''-f(u) = g , \quad u(0) = u(\pi) = 0   \ .
\end{equation}

The nonlinearities we have in mind  may be asymptotically or piecewise linear. The function $ f_{a\ell}: \RR \to \RR$ is {\it asymptotically linear} if it is a $C^1$, convex function which is asymptotically linear with parameters $\ell_-$,  $\ell_+$,
\[
\lim_{x \to - \infty}  f_{a\ell}(x) = \ell_-,  \quad \lim_{x \to \infty}  f_{a\ell}(x) = \ell_+ \ ,\]
On the other hand,  $f_{p\ell}$ is {\it piecewise linear} (\cite{COSTA, LAZERANDMCKENNA, TELES}) if it is of the form
\begin{equation} \label{piecewise}
f_{p\ell}(x) = \begin{cases}
	\ell_- x\ , \ x < 0   \\
	\ell_+ x\ , \ x > 0
\end{cases} \ .
\end{equation}

\subsection{The continous case:  $F_{a\ell}$ and $F_{p\ell}$} \label{contSturm}

We  consider  operators $F_{a\ell}(u) = - u'' -  f_{a\ell}(u)$ and $F_{p\ell}(u) = - u'' - f_{p\ell}(u)$. Clearly, $F_{p\ell}$ is $1$-homogeneous: $F_{p\ell}( t u ) = t F_{p\ell}(u)$ for $t >  0$.

For the linearities we have in mind, there are two natural choices for domain $X$ and counterdomain $Y$ of $F: X \to Y$. For $X$, one may take the subspaces of $C^2([0,\pi])$ or $H^2([0,\pi])$ of functions satisfying Dirichlet boundary conditions. Correspondingly, $Y$ should be either $C^0([0,1])$ or $L^2([0,1])$. The results in this section apply to both choices of pairs $X$ and $Y$.

The critical set of the operator $F_{a\ell}$   has been studied in \cite{BURGHELEA} and \cite{TOMEIEBUENO}. Their Jacobians are self-adjoint operators with simple spectrum, which are bounded from below:  we label eigenvalues in increasing order.

\medskip

Recall that the linear operator $u \mapsto - u''$ acting on functions satisfying Dirichlet boundary conditions has eigenvalues equal to $\lambda_k = k^2, k=1, 2, \ldots$ and associated eigenvectors $\phi_k = \sin(kx)$. For special functions $g$,
one may count solutions of $F(u) = g$, as follows.

\medskip
\begin{theorem}[Costa-Figueiredo-Srikanth \cite{COSTA}] \label{Fp}
For parameters $\ell_-$,  $\ell_+$
satisfying
\[
\ell_- < \lambda_1 = 1, \ \quad \lambda_k = k^2<\ell_+<(k+1)^2 = \lambda_{k+1} \ . \]
Then,  the equations $F_{p\ell}(u) = -\sin(x) $ and $F_{a\ell}(u) = -t\sin(x) $ for large $t > 0$
have exactly $2k$ solutions.
\end{theorem}

\bigskip

\subsection{The discrete case:  $F^h_{a\ell}$ and $F^h_{p\ell}$} \label{discSturm}
\medskip
We  discretize $F_{a\ell}$ and $F_{p\ell}$
on the regular mesh
\[ \overline{I_h}=\left\{x_i = ih \ , \ i=0,...,n+1 \ ,   h=\frac{\pi}{n+1} \right\} \ .\]
For interior points in $I_h=\overline{I_h} \setminus \{x_0,x_{n+1}\}$,  we take the usual approximation
$$
-u''(x_i) \sim \frac{1}{h^2} \left(-u(x_{i+1})+2u(x_i)-u(x_{i-1})\right)
$$
and define the tridiagonal,  $n \times n$ symmetric  matrix $A^h = (1/h^2) A$, where $A_{ii} = 2, i=1...n$ and $A_{i,i+1}= A_{i+1,i} = -1, i=1...n-1$. The  eigenvalues of $A^h$ are   \[ \lambda_k^h \ = \frac{2}{h^2}( 1 - \cos(kh))  , \ k = 1, \ldots, n \]
with associated eigenvectors $\phi_k^h = \sin(k I_h)= (\sin(kx_i))$ for $x_i \in I_h$.
Similarly, set $u^h = u(I_h), f(u^h)=(f(u(I_h))^T \in \RR^n$. The discretized operators act on $\RR^n$,  
\[ F_{a\ell}^h(u^h) = A^h u^h -  f_{a\ell}(u^h) \ , \quad F_{p\ell}^h(u^h) =  A^h u^h - f_{p\ell}(u^h) . \]

Here is the discrete counterpart of Theorem \ref{Fp}. A vector $p \in \RR^n$ is {\it positive},  $ p >0$, if all its entries are strictly positive. 

\medskip

\begin{theorem}[Teles-Tomei \cite{TELES}] \label{contando} Let $F^h_{a\ell}$ or $F^h_{p\ell}$ as above, with parameters $\ell_-$ and $\ell_+$ satisfying
\[  \ell_-< \lambda_1^h   \ ,  \quad \lambda_k^h < \ell_+  < \lambda_{k+1}^h \ . \] 
Let $y,p \in \mathbb{R}^n $, $p > 0$. Then, for fixed $  \ell_-< \lambda_1^h$, large  $t >0$ and large $\ell_+ $, the equations
\[F_{a\ell}^h(u^h)=y -tp \ \ \hbox{and} \ \ F_{p\ell}^h(u^h)=y -tp ,\]
have exactly $2^n$ solutions.
\end{theorem}

\bigskip

As an example, let $n = 15$, $h = \pi/16$,  $g  = - \sin(I_h)$ and $I_h$ as above. The extreme eigenvalues of $ A^h$ are $\lambda_1  = 0.99679136$ and $\lambda_{15}  =102.7561006$.  Consider different choices of $f$ as in (\ref{piecewise}) with  parameters
\[ \ell_-=\frac{\lambda_{1} }{2} \ \ \text{and} \ \  \ell_{+}^k=\frac{\lambda_k +\lambda_{k+1} }{2}, \ \text{for} \ k=1,...,14 \ ,  \ \ \ell_{+}^{15}=\lambda_{15} +\frac{\lambda_{1} }{2}.\]
The number of solutions $N$ of $ F  (u ) = g $ for different $\ell_+^k$ is given below (\cite{TELES}).

\begin{table} [ht] \label{tabela}  
\small
\centering
\begin{tabular}{|c|c|c|c|c|c|c|c|c|c|c|c|c|c|c|c|}
	\hline
	$k$ &1&2&3&4&5&6&7&8&9&10&11&12&13&14&15  \\ \hline
	$N$ & 2 & 4 & 6 & 8 & 12 & 12 &22 & 24 & 32 & 100 & 286 &634 & 972 & 1320 & 2058\\ \hline
\end{tabular}
\label{Ng}
\end{table}

For $\ell_{+} = \lambda_{15}  + 100 \lambda_1  \sim 202.4352$, there are $2^{15}$ solutions, in accordance with Theorem \ref{contando}. 
For small $k$, $N = 2k$, as for $ F_{a\ell}$ in the previous section. As $k$ increases, $N$ becomes exponential.

\medskip

The abundance of solutions of $F^h_{p\ell}(u^h) = g^h$ makes it an interesting test case for our strategy. Most such solutions  do not admit a continuous limit.

\medskip

For our algorithm, this is the only required knowledge about the critical set: 

\medskip
\noindent{\bf Main Property:} Straight lines in the domain which are parallel to a positive function (a positive vector, in the discrete case) intercept the critical set  abundantly. 

\smallskip

The property is proved in a more general setting (for Jacobians of the form $DF(u) = - \Delta + f'(u)$) in  Proposition   \ref{verticais}.

\bigskip
How can there be so many solutions? The geometry of $\cC$ and $F (\cC)$ of the operators $F  = F^h_{a\ell}$  illuminates the situation. Again, this information is not relevant for application of the strategy. We synthesize \cite{TELES}.

Given the smallest eigenvalue $\lambda_1 $ of $A $ and its (normalized, positive) eigenvector $\phi_1 $, set $V = \langle \phi_1  \rangle$ and split  $\mathbb{R}^n = V^\perp \oplus  V$. For $\ell_+ >> 0$, the critical set $\mathcal{C}$ of $F : \mathbb{R}^n \to \mathbb{R}^n$ consists of $n$ hypersurfaces $\mathcal{C}_j, j=1, \ldots,n$.  
Each $\mathcal{C}_j$ is a graph of a function $c_j : V^\perp \to V$ which projects diffeomorphically from $\mathbb{R}^n$ to $V^\perp$. Thus, topologically, $\mathcal{C}_j$ is trivial and 
each line $r$ parallel to $\phi_1 $ intercepts all the critical components $\mathcal{C}_j$ of $F $.

The images $F(\mathcal{C}_j)$ are more complicated. Suppose the interval $(\ell_-,\ell_+)$  contains all eigenvalues $\lambda_1< \ldots < \lambda_n$. Then $F(\mathcal{C}_j)$ wraps around the line
$\{t(1, 1, \ldots, 1), t \in \mathbb{R} \}\subset \mathbb{R}^n$ substantially: $\binom{n-1}{j-1}$ times! It is this  `topological turbulence' which gives rise to the abundance of preimages for appropriate right hand sides. Indeed, from item (4) of Theorem \ref{theo:model}, points at opposite sides close to $F(\cC)$  have their number of preimages differing by two. This situation should compared  with the second function in Section \ref{Visual}: cubic behavior at infinity leads to winding of the image of the outermost critical curve and to nine zeros of $F$.

\bigskip

\subsection{Piecewise linear geometry: the critical set of $ F^h_{p\ell}$} \label{table}

In this subsection, we write $F$ to denote $F^h_{p\ell}$, $f$ for $f_{p\ell}$ and drop the discretizing parameter $h$
on $A, \phi_1, u$ and $g$. Now, for $t > 0$, $F(t u) = t F(u)$.


As an application of our strategy, we find solutions of
$ F_p (u ) = g = - \sin(I_h)$,
for $n=15$ (so that $h = \pi/16$),   and $f$ with different choices of parameters $\ell_-$ and $\ell_+$. 
Given $ F $ and $g $, the solutions of $ F (u ) = g $ are obtained explicitely rather trivially. Here, we compare these answers with those generated by our strategy.

\bigskip

As  vector coordinates in each (open) orthant $\cO$ have fixed signs, the restriction of the function $ F $ to  $\cO \subset \RR^n$ is linear. Its extension to $\RR^n$ is clearly continuous (but not differentiable at coordinate hyperplanes, the subspaces in the boundary of two quadrants). The numerics in \cite{TELES} generates a reliable data bank of solutions: solve a linear system in each orthant and check if the solution  belongs to it.

We now consider the critical set $\cC$ of $F = F^h_{p\ell}$.  On an (open) orthant $\cO \subset \RR^n$, the function $f$ is defined by the sign of the entries of its vectors: for $v \in \cO$, $f(v)$ is multiplication by a diagonal matrix $D^\cO$, with diagonal entries $d_i$ equal to $\ell_-$ or $\ell_+$, according to the sign of entry  $v_i$. Thus, the continuous function $F : \RR^n  \to \RR^n$ is of the form $F  = A  - D^{\cO}$ when restricted to $\cO$.

Generically (i.e., for an open, dense set of pairs $(\ell_-, \ell_+)$),  the matrices $F = A  - D^{\cO}$ are invertible for all $\cO$. If this is not the case, full orthants belong to the critical set $\cC$ and then $\cC$ has nonempty interior.  Assuming invertibility,  $F$ is trivially a local diffeomorphism in each orthant $\cO$. The critical set of $F $ makes (topological) sense: it is the set in which $F $ is not a local homeomorphism.

Say a vector $v$ is in the boundary of exactly two orthants, $\cO_1$ and $\cO_2$ (equivalently, $v$ has a single entry equal to zero). If $\sign \det (A  - D^{\cO_1}) = \sign \det (A  - D^{\cO_1})$, the function $F $  is a local homeomorphism at $v$. If instead $\sign \det (A  - D^{\cO_1}) =  -\sign \det (A  - D^{\cO_1})$, and $F $ near $v$ is a topological fold. 
A {\it slab} is the closure of the subset of a coordinate hyperplane in the boundary of two quadrants whose vectors have signs at exactly one coordinate.
The critical set consists of slabs across which $\det (A  - D^{\cO})$ changes sign. In the nongeneric situation, whenever $A  - D^{\cO}$ is not invertible, we  include $\cO$ in the critical set.

\bigskip

The case $n=2$ is visualizable. The  matrix $A $
has eigenvalues  $\lambda_1 \approx 0.9119$ and $\lambda_2 \approx 2.7357$. For two pairs $(\ell_-, \ell_+) = (-1, 2)$ and $(-1, 4)$, 
Figure 7 
describes the critical sets (slabs are half-lines through the origin) and their images. In each quadrant of the domain we indicate $\sign \det (A  - D^\cO)$. In each connected component of  $\mathbb{R}^2 \setminus F (\mathcal{C}) $ we specify instead the number of preimages, in accordance with item (4) of Theorem \ref{theo:model}. 

The fact that, for $\ell_+ = 4$, $\cC$  is not a (topological) manifold, is innocuous. Indeed, for the general operator $F$, the critical set $\cC$ is connected, as slabs always contain the origin. In particular, the critical set is rarely a (topological) manifold: at the origin, many pieces of the coordinate hyperplanes may meet.

\begin{figure}  [ht] 
\centering
\begin{minipage}{0.45\textwidth}
	\centering
	\includegraphics[width=1.0\textwidth]{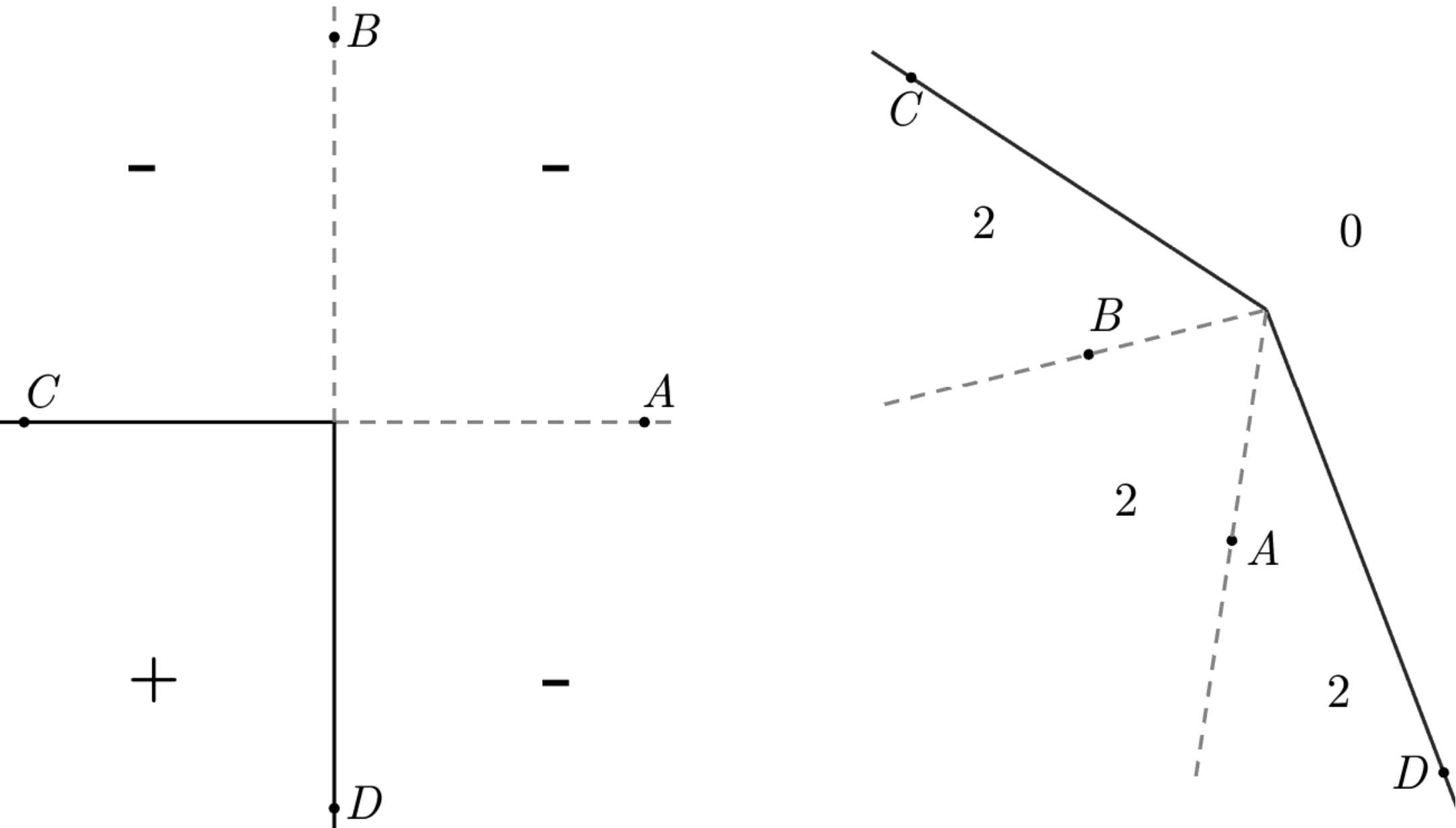} 
	\label{fig4a}
\end{minipage}
\begin{minipage}{0.45\textwidth}
	\centering
	\includegraphics[width=1.0\textwidth]{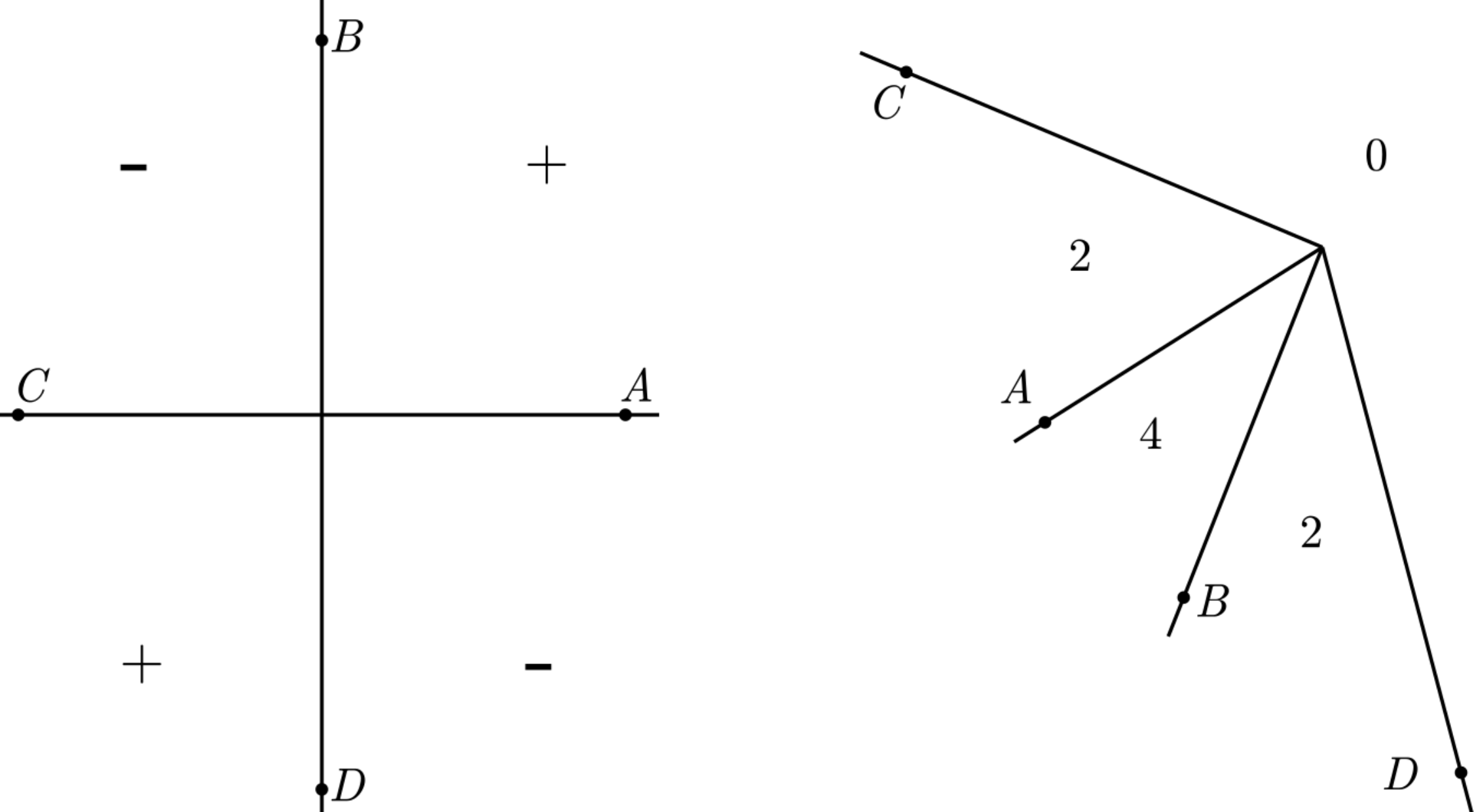} 
\end{minipage}
\hfill

\caption{$\mathcal{C}$ and $F (\mathcal{C})$ for  $\ell_-=-1, \ell_+=2$ and for $\ell_-=-1, \ell_+=4$}
\end{figure}

Back to the general case ($n$ arbitrary), as eigenvalues increase along a (generic) line $r$ parallel to the vector $\phi_1 $ (a consequence of  Proposition  \ref{minmax}), we must have that $r$  intercepts $\mathcal{C}$ in as many eigenvalues of $A $ there are between $\ell_-$ and $\ell_+$.
Notice the affinity between the lines $r$ and the rays from the origin in the examples of Section \ref{Visual}: as in Section \ref{winding}, different lines in the domain may contain a different set of solutions in their bifurcation diagrams.

We consider again $n=2$ with $(\ell_-, \ell_+) = (-1, 4)$. Imitating the continuous case $F_{p\ell}$ in \cite{LAZERANDMCKENNA},   two solutions of $ F  u  = F^h_{p\ell}  u  =  A  u  - f(u )= -t\phi_1 (x)$ are
\begin{equation} \label{LazerMcKenna}
\frac{t\ \phi_1 }{\ell_+ -\lambda_1 } > 0 \ \ \
\text{and} \ \ \
\frac{t\ \phi_1 }{\ell_- - \lambda_1 } < 0,  \quad \lambda_1 =\frac{2}{h^2} \left(1-\cos h \right) \quad \hbox{for}\ \  \ell_- < \lambda_1  < \ell_+ \ .
\end{equation}
Through the positive solution
$P_0 \approx ( 280.4396, 280.4396)$ obtained from \eqref{LazerMcKenna} by setting $\ell_+ = 4$, draw the half-line
$r =\{P_0+s(0.2\sin(2I_h)-0.8\sin(I_h)), \  s \geq 0\} \subset \mathbb{R}^{2}$.
We add the term  $0.2\sin(2I_h)$ to minimize the possibility that $r$ intercepts $\cC$ at points which are not (topological) folds  (i.e., points with more than one entry equal to zero -- this is especially relevant for $n$ large). Figure 8 
shows $r$, $F (r)$ and the bifurcation diagram $\mathcal{B}$, the connected component of $F^{-1}(F(r))$ containing $P_0=(0,0)$.
The critical points  in $r$ are $a_*$ and $c_*$. At these points, $F$ is a fold, and  $r$ is mirrored across $\cC$.  These reflections (bifurcation branches) in turn  intercept $\cC$ (at $b_*$). By continuation, the missing three solutions $P_i, i=1,2,3,$ are obtained.

\begin{figure}  [ht] \label{fig4b}
\centering
\begin{minipage}{0.45\textwidth}
	\centering
	\includegraphics[width=0.9\textwidth]{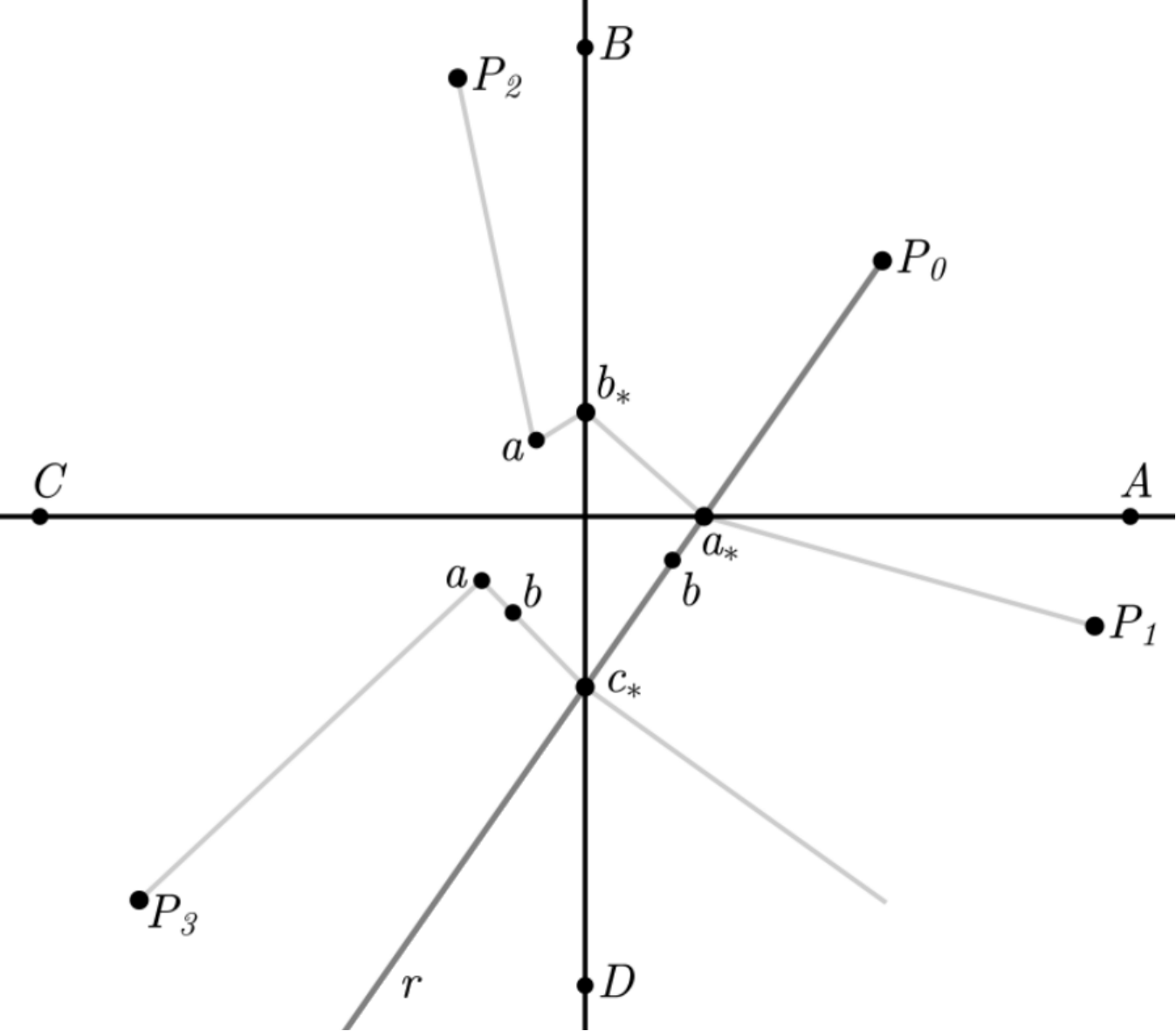} 
\end{minipage}
\begin{minipage}{0.45\textwidth}
	\centering
	\includegraphics[width=0.75\textwidth]{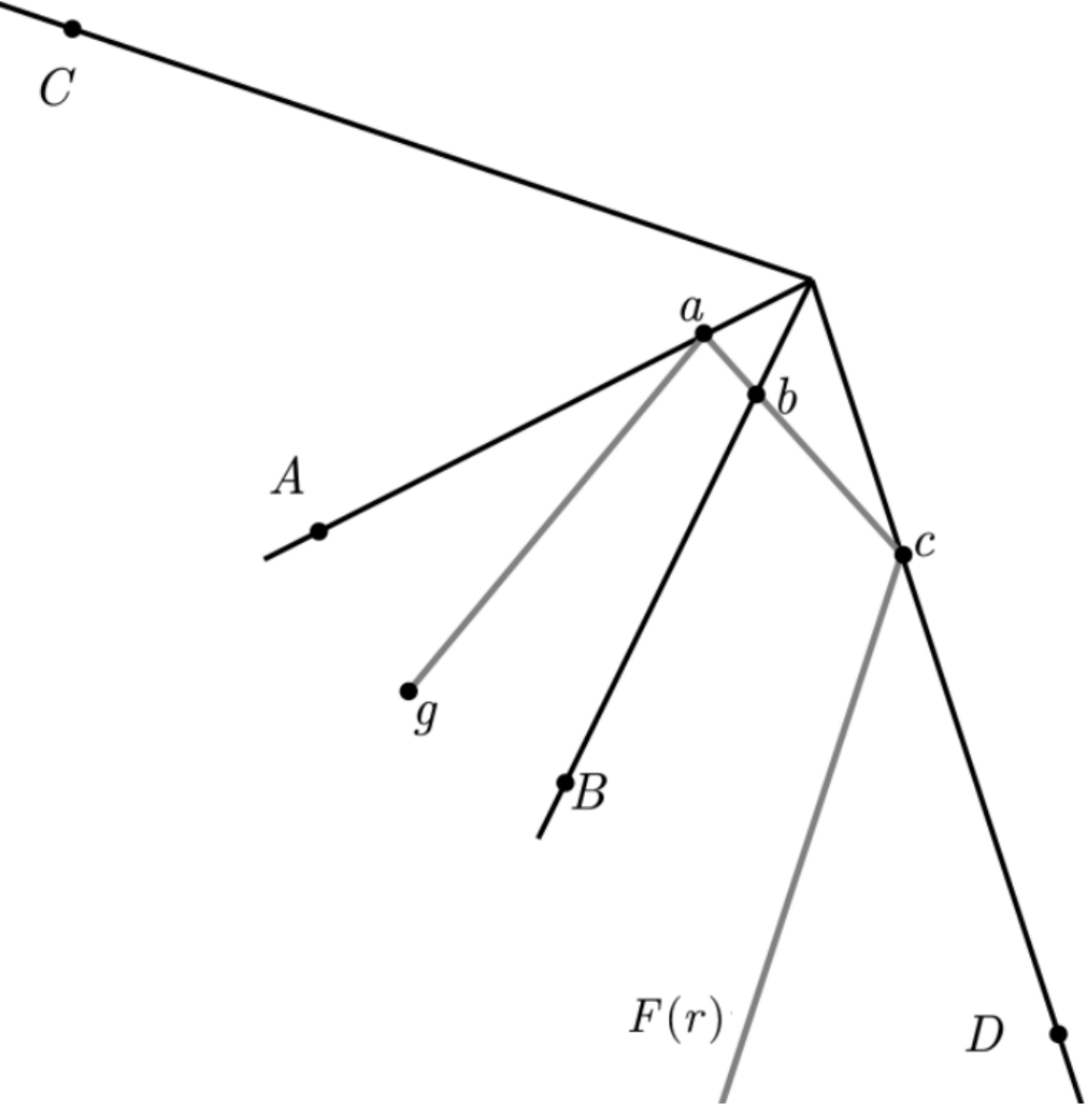} 
\end{minipage}
\caption{$\mathcal{B}=F^{-1}(F(r))$ (left) and $F (r)$ (right)  for $\ell_-=-1, \ell_+=4$}
\hfill
\end{figure}

%

We finally consider the case $n=15$ for two choices of asymptotic parameters.

\begin{figure}  [ht]
\centering
\begin{minipage}{0.45\textwidth}
	\centering
	\includegraphics[width=165pt]{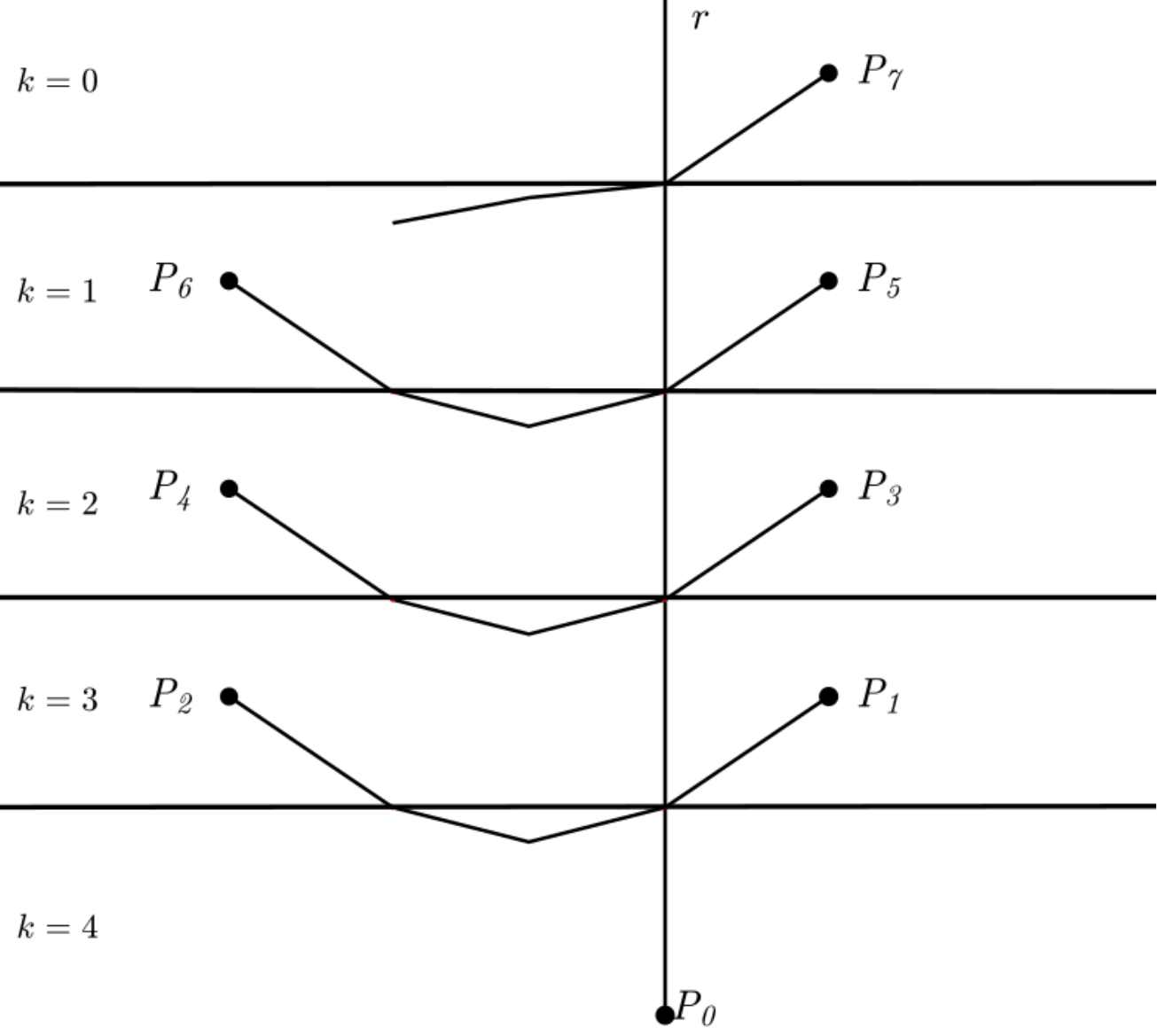}
	\caption{$\mathcal{B}$ for $k=4$.}
	\label{fig:4.7}
\end{minipage}
\hspace{0.5cm}
\begin{minipage}{0.45\textwidth}
	\centering
	\includegraphics[width=165pt]{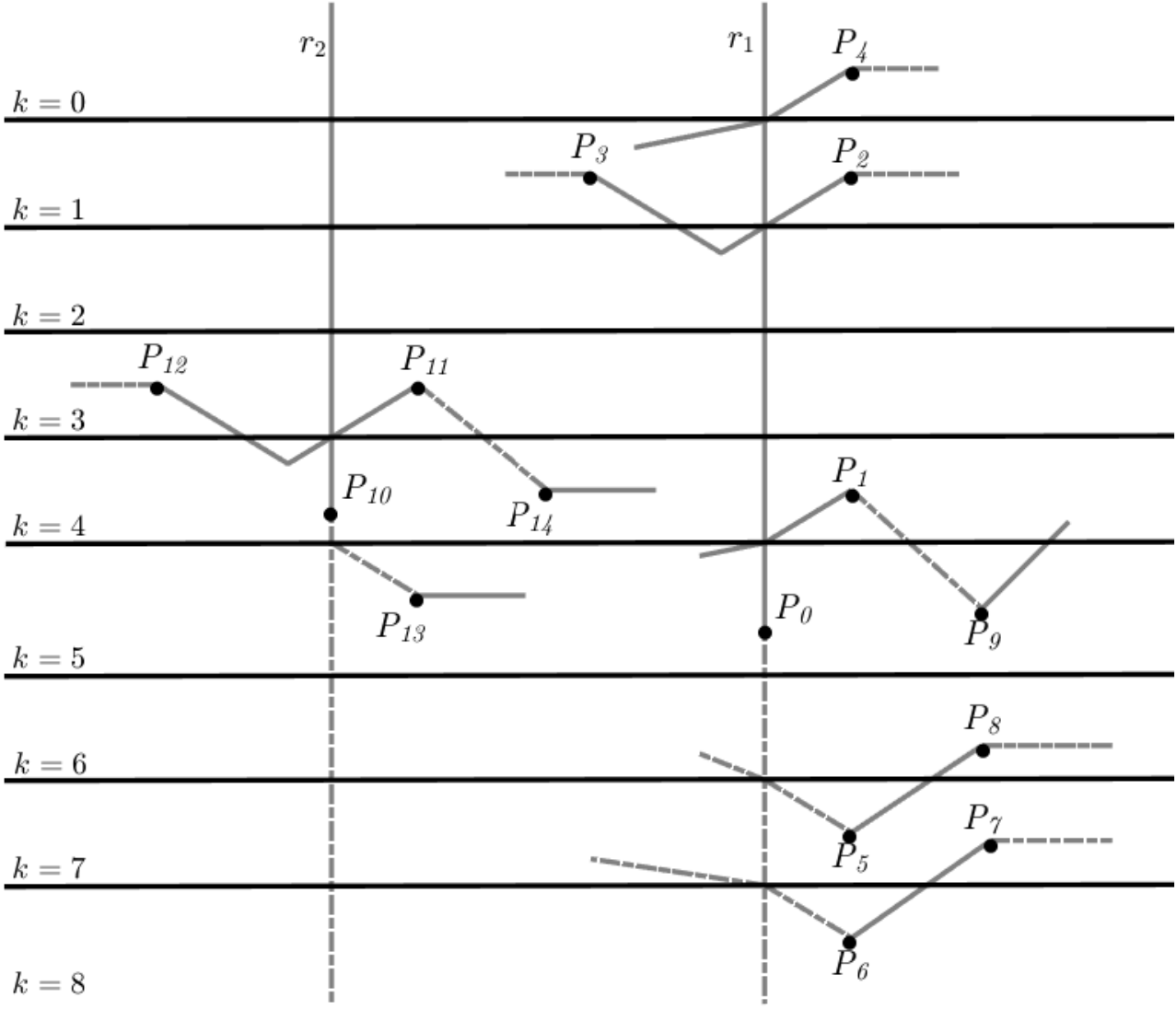}
	\caption{Two $\mathcal{B}$'s for $k=8$}
	\label{fig:twolines}
\end{minipage}
\hfill
\end{figure}

The pair $(\ell_-, \ell_+) = ( 0.4984 ,19.1248)$ encloses four eigenvalues of $A $ (in Table \ref{tabela}, $k=4$). 
A positive solution of $F (u) = g =-\sin(I_h)$ is $P_0=0.0551633\sin(I_h)$, from \eqref{LazerMcKenna}. We use half-lines which are perturbations of the eigenvector $\sin(I_h)$ associated with the smallest eigenvalue $\lambda_1 $,
\[ r=\{P_0 + s(\sin(I_h)-0.1\sin(2I_h) -0.1\sin(3I_h)- 0.1\sin(4I_h)), \ \text{com} \ s \geq 0\} \subset \mathbb{R}^{15} ,\] 
to ensure transversal intersections of $r$ with the critical set $\mathcal{C}$ of $F $. Figure  \ref{fig:4.7} represents schematically  the five open subsets  $R_k, k= 0, \ldots,4$ of the regular set, $\RR^{15} \setminus \mathcal{C}$, each containg the Jacobians with  $k$ negative eigenvalues. Lines separating the regions are critical points.
The bifurcation diagram $\mathcal{B}$ associated with $r$ from $P_0$  contains all the preimages in the data bank associated with Table \ref{tabela}.

\medskip


\medskip

We now take  parameters $(\ell_-, \ell_+) = (0.4984,56.9367)$ and use the same $g$. There should be  24 preimages from Table \ref{table} (here, $k=8$). We initialize our strategy differently: choose orthants randomly, and search for a first solution by solving the associated linear system. A first solution $P_0$ (not necessarily given by equation (\ref{LazerMcKenna})), was obtained after drawing 127 orthants (out of $2^{15}$). Consider the (possibly generic) line through $P_0$,
\begin{align}
r_1=&\{ P_0+ s(0.8\sin(I_h)-0.1\sin(2I_h)-0.1\sin(3I_h)-0.1\sin(4I_h)\nonumber \\ &-0.1\sin(5I_h)-0.1\sin(6I_h)-0.1\sin(7I_h)+0.1\sin(8I_h))\} .\nonumber
\end{align}

The  bifurcation diagram $\mathcal{B}$, shown partially on the right of Figure \ref{fig:twolines}, led to 17 new candidate preimages, of which only 9 had a small relative error, of the order of $10^{-15}$. The spurious solutions are related to inversion
of unstable matrices of the form $(A -D^o)$ in orthants trespassed along the homotopy process. Stretches of the bifurcation diagram are  continuous or dotted, depending if   $s$ is positive or negative. 

A new solution $P_{10}$ is obtained by sampling additional orthants. The bifurcation diagram associated with $r_2=\{P_{10}+s(\sin(I_h)-\sin(8I_h))\}$,  on the left of Figure \ref{fig:twolines}, yielded four more preimages, after filtering candidates by relative error. 


Additional solutions were obtained by sampling, leading to $24$ solutions. The table below counts solutions $u$ by  the number of negative eigenvalues of $DF(u)$.

\begin{table} [h]
\centering
\begin{tabular}{|c|c|c|c|c|c|c|c|c|c|}
	\hline
	$k$ & 0 & 1 & 2 & 3 & 4 & 5 & 6 & 7 & 8 \\ \hline
	$\hbox{Number of preimages}$ & 1 & 2 & 2 & 4 & 6 & 4 & 2 & 2 & 1 \\ \hline
\end{tabular}
\end{table}

%


\section{Semi-linear perturbations of the Laplacian} \label{Solimini}

To introduce the basic ideas, we present the simplest nontrivial flower $\cF = F^{-1}(F(\cC))$ associated with a differential operator $F$ with critical set $\cC$.
The following version of the celebrated Ambrosetti-Prodi theorem combines material from \cite{AP,MM, BERGER, CALNETO, SMILEY}).  On a bounded set $\Omega \subset \mathbb{R}^n$ with smooth boundary, let  
\[ -\Delta_D: X = H^2(\Omega) \cap H^1_D(\Omega) \to Y=L^2(\Omega)\]
be the Dirichlet Laplacian, and denote its smallest eigenvalues by $\lambda_1$. A {\it vertical line} is a line in $X$ or $Y $ with direction given by $\phi_1>0$, a positive eigenvector associated with $\lambda_1$. {\it Horizontal subspaces} $H_X \subset X $ and $H_Y\subset Y$  consist of vectors perpendicular to $\phi_1$. {\it Horizontal  hyperplanes}  are parallel to the horizontal subspaces.

\medskip
\begin{theorem} 	Consider the function
\begin{equation*} \label{AP}
	F: X \to Y \ , \quad u \mapsto  - \Delta_D u - f(u) 
\end{equation*}
where $f: \RR \to \RR$ is a  strictly convex smooth function satisfying
\begin{equation} \label{ASY} - \infty < \lim_{x \to - \infty} f'(x) < \lambda_1 < \lim_{x \to  \infty} f'(x) < \lambda_2 \ .
\end{equation}
Then the function $F$ is proper and its critical set $\cC$ contains only folds.
The orthogonal projection $\cC \to H_X$ is a diffeomorphism, as is the projection of the image of each horizontal hyperplane $F(H_X + x) \to H_Y, x \in X$. 
The inverse under $F$ of vertical lines in $Y$ are (connected) curves in $X$  intercepting each horizontal hyperplane and $\cC$ exactly once, transversally. In particular, $F$ is a global fold and the equation $F(u) = g \in Y$ has 0, 1 or 2 solutions.
\end{theorem}

\medskip
From Section \ref{basicfolds},  $F: X \to Y$ being a {\it global fold} implies the existence of
diffeomorphisms $\Phi: \RR \times Z \to X $ and  $\Psi: Y \to \RR \times Z$ such that
\[ \tilde F = \Psi \circ F \circ \Phi(t, z): \RR \times Z \to \RR \times Z  \ , \quad \tilde F (t, z) =  (t^2, z)\] for some real Banach space $Z$.  In particular, $\cC$ and $F(\cC)$ are diffeomorphic to hyperplanes, and  $F^{-1}(F(\cC)) = \cC$. Compared to the examples in the previous sections, the global geometry of $F$ is very simple.  Both functions $F$ and $\tilde F$ trivially satisfy Theorem \ref{theo:model}: domain and codomain split in two tiles, both topological half-spaces, and both domain tiles  are sent to the same image tile, as in Figure	\ref{fig:extensao2}. In both cases, the flower coincides with $\cC$.

The simple geometry led to numerical approaches to solving $F(u) = g$.
Smiley (\cite{SMILEY}) suggested an algorithm based on one-dimensional searches, later implemented in \cite{CALNETO}. A similar finite dimensional reduction applies for (generic) asymptotically linear functions $f$  for which the image of $f'$ contains a finite number of eigenvalues of $-\Delta_D$ (\cite{CALNETO,KAMINSKI}).

Under hypothesis (\ref{ASY}), the {\it nonconvexity} of $f$ implies that some right hand side $g$ has four preimages (\cite{CaTZ}). Up to technicalities, our strategy yields {\it all} solutions for convex and nonconvex nonlinearities.

\medskip
Breuer, McKenna and Plum (\cite{PLUM}) computed four solutions  of
\begin{equation*} \label{Plum} F(u) = - \Delta u + u^2 \ = \  800 \sin( \pi x) \sin(\pi y), \ (x,y ) \in \Omega = (0,1) \times (0,1) \ , \ u|_\Omega = 0  
\end{equation*}
making use of a variational formulation of the equation. One solution, call it $u_s$, is  a saddle point  of a functional $\Phi$ associated with $F$.  The authors present a computer assisted proof that $u_s$ is reachable by Newton's method from a computed initial condition $\tilde u_s$.  Such four solutions were also obtained in \cite{ALLGOWER2}. The equation  may be treated with our methods, but we introduce a situation with additional difficulties.

\bigskip
In our final example, we consider $F: X \to Y, F(u) = - \Delta_D - f(u)$ with a special nonlinearity $f$. Again, $\Omega \subset \RR^n$ is a bounded, smooth domain, let $X = H^2(\Omega) \cap H^1_0(\Omega), Y = L^2(\Omega)$.
As in  previous sections, we count intersections of  vertical lines and $\cC$, the critical set of $F$ ( Proposition \ref{verticais}).

\medskip
The result below is the standard   min-max theorem (Theorem XIII.1, \cite{RS}) for a self-adjoint operator $T: X \to Y$ which is bounded from below, containing only point spectrum. Label the eigenvalues $\lambda_i = \lambda_i(T)$ in nonincreasing order, counting multiplicity, with associated eigenvectors $\phi_i = \phi_i(T)$.

\medskip

\begin{proposition} \label{minmax} Set $S = \{ v \in X \ | \ \langle v , v \rangle  = 1\}$. Let  $X_{i-1} \subset X$ a subspace generated by $i-1$ linearly independent eigenvectors $\phi_1, \ldots, \phi_{i-1}$ and denote by $X_{i-1}^\perp \subset X$ its ($L^2$) orthogonal complement.  Then, for $T: X\to Y$ and $\lambda_i$ as above, we have
\[ \lambda_i = \sup_{X_{i-1}} \ \inf_{v \in X_{i-1}^\perp \cap S} \langle v, T v\rangle \ .\]
Moreover, for each $i \in \NN$, there is a  bounded set $B_i \subset X$ (in the $X$-norm) for which 
\[\lambda_i = \sup_{X_{i-1}} \  \inf_{v \in X_{i-1}^\perp \cap S \cap B_i} \langle v, T v\rangle \]   
\end{proposition}

\begin{proof}
The first expression for $\lambda_i$ is familiar. For the second claim, 
let $X_i$ be the subspace generated by the span of $X_{i-1}$ and $\phi_i \in X_{i-1}^\perp$, a normalized eigenvector associated with $\lambda_i$. The intersection $B_i = X_i \cap S$ is a bounded subset of $X_i$, both in the $L^2$ and $X$-norm, containing $\phi_i$. 
\end{proof}

For $w \in X$, consider the {\it vertical line}  $r(t)= \{w + t \phi_1,  t \in \RR\} \subset X$. Let $\overline{\RR} = \RR \cup \{ -\infty, \infty\}$. For  $ t \in \overline{\RR}$, define $T[t]: X \to Y$ as
\[ T[-\infty] v = - \Delta_D v - \ell_- v ,  \quad T[t] v \ = \   - \Delta_D v - f'(r(t)) v , \quad   T[\infty] v = - \Delta_D v - \ell_+ v . \]
The next proposition is essentially in Appendix 5.2, \cite{CaTZ}.

\medskip
\begin{proposition} \label{continuity} For each $i \in \NN$, the functions 
$ \lambda_i:  \overline{\RR} \to \RR ,  t \mapsto \lambda_i(T[t])$ 
are continuous. If  $f$ is a strictly convex function, then  the  $\lambda_i$'s are strictly decreasing.
\end{proposition}

\begin{proof} Fix $i \in \NN$ and let $t \to t_0 \in \overline{\RR}$.
We abuse notation and set $\ell_\pm = f'(r(t_0))$ for $t_0 = \pm \infty$.
From Proposition \ref{minmax}, it suffices to show  $\langle T[t] v, v \rangle \to \langle T[t_0] v, v \rangle$, or equivalently,  $\langle f'(r(t)) v, v \rangle \to \langle f'(r(t_0)) v, v \rangle$, uniformly in $v \in X \cap S \cap B_i$, for some given bounded set $B_i \subset X$.

Now, $X$ embeds continuously on $L^{2n/(n-4)} \subset L^4$ for $n \ge 5$  and in every $L^p$ for $n < 5$. We then have, for $\rho= (4/3), \sigma=4$, 
\begin{align*} 
	| &\langle (f'(r(t))- f'(r(t_0)) v, v \rangle| \  \le \ \int_\Omega |(f'(r(t))- f'(r(t_0)) v^2 |\ dx \\ &\le \|f'(r(t))- f'(r(t_0)) \|_{L^\rho} \  \| v \|_{L^{2\sigma}}  \le C \|f'(r(t))- f'(r(t_0)) \|_{L^\rho} \  \| v \|_{H^{2}} \ .
\end{align*}
Since $v \in B_i$, $ \| v \|_{H^{2}}$ is uniformly bounded, and it is easy to see that, as $t \to t_0$ we have $\|f'(r(t))- f'(r(t_0)) \|_{L^\rho} \to 0$. This proves continuity. 

Strict monotonicity arises from the fact that critical points of the quadratic form in Proposition \ref{minmax} are eigenfunctions of $T[t]$. The finite multiplicity of each $\lambda_*$ provides compactness of the set of normalized associated eigenvectors which in turn guarantees that all eigenvalues equal to $\lambda_*$ increase.
\end{proof}

We do not have  $\|f'(r(t))- f'(r(t_0)) \|_{L^\infty} \to 0$:  for $t_0 =  \infty$, on  $\partial \Omega$, $f'=0 \ne \ell_+$.

\medskip
\begin{proposition} \label{verticais} Let $f: \RR \to \RR$ be a smooth,  function for which 
\begin{equation*} \label{ASYk} - \infty < \ell_- = \lim_{x \to - \infty} f'(x) < \lambda_1 < \ldots \le \lambda_k <  \ell_+ = \lim_{x \to  \infty} f'(x) < \lambda_{k+1} \ .
\end{equation*}
Then vertical lines meet the critical set $\cC$ of $ F$  at at least $k$ points (counted with multiplicity). If $f$ is strictly convex, there are exactly $k$ intersection points. 
\end{proposition}

\begin{proof}
The eigenvalues of $T[-\infty]$ and $T[\infty]$ are respectively 
\begin{equation*} \label{extremos}
	\begin{aligned} &0 < \lambda_1^D - \ell_1 \le \lambda_2^D - \ell_1 \le \ldots \le \lambda_k^D - \ell_1 <  \lambda_{k+1}^D - \ell_1 \le \ldots \to \infty \ , \\ &\lambda_1^D - \ell_k \le \lambda_2^D - \ell_k \le \ldots \le \lambda_k^D- \ell_k < 0 <  \lambda_{k+1}^D- \ell_k \le \ldots \to \infty \ . 	\end{aligned}
\end{equation*}
Fix a vertical line $r = \{w + t \phi \}$. 
From Proposition \ref{continuity},  each $\lambda_i(T[t])$ along $r(t)$ is continuous in $\overline{\RR}$. Thus, as $t$ varies between extremal values, each of the first $k$ eigenvalues must be zero at (at least some) some point in $X$. Such intersection points are exactly the critical points of $F$ in the vertical line $\ell$, as $\cC$ consists of the points for which some eigenvalue of the Jacobian $DF(u) v = - \Delta_D v f'r(t))v$ is zero.
\end{proof}

\begin{figure} [ht]

\begin{centering}
	\includegraphics[height=150pt,width=150pt]{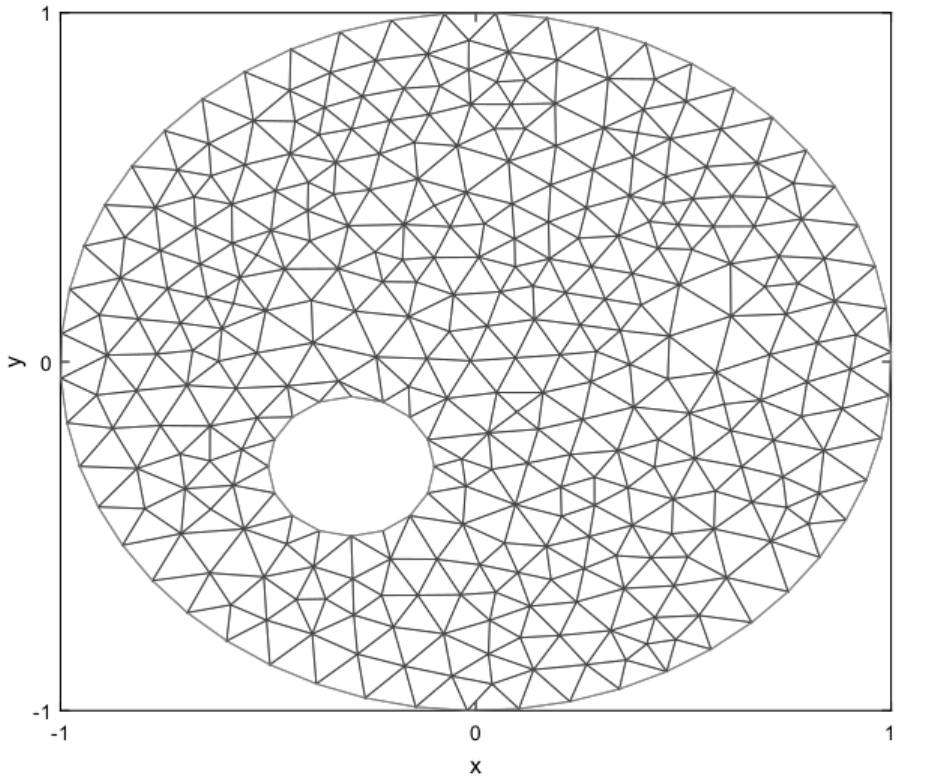}
	\caption{The domain $\Omega$ and its mesh}
	\label{fig:5.1}
\end{centering}
\end{figure}

We are ready to present our last example.
In \cite{LMCKENNA}, Lazer and McKenna conjectured that, for an asymptotically linear $f$ with parameters $\ell_-<\lambda_1$ and $\lambda_k<\ell_+<\lambda_{k+1}, k \in \mathbb{N}$, there should be at least  $2k$ solutions of (\ref{eq:PVCnaolinear}) for  $g = - t \phi_1, t >> 0$. A counterexample was provided by Dancer
(\cite{LAZERANDMCKENNA}). The example we consider, from \cite{DIEGO}, follows a positive result of Solimini \cite{SOLIMINI}. 

Let  $\Omega$ be the annulus \[ \Omega=\{x \in  \mathbb{R}^2 | \ |x|<1 \ ,  \ |x-(-0.3,-0.3)|>0.2 \} \ . \] 

We discretize  $- \Delta_D$
by piecewise linear finite elements on a mesh on $\Omega$ with 274 triangles, as in Figure \ref{fig:5.1}.
The four smallest eigenvalues of the discretized operator $-\Delta_D $ are simple,
\[ \lambda_1  \approx  9.0988, \quad \lambda_2  \approx 16.3218, \quad \lambda_3  \approx 22.9346, \quad \lambda_4  \approx 30.4949 \ . \]

We consider
\begin{equation} 	\label{eq:PVCnaolinear}
F(u)=-\Delta_D u-f(u)=g, \ \ \  \ u|_{\partial \Omega}=0 .
\end{equation}

According to Solimini, for some $\epsilon>0$ and parameters $\ell_-$ and $ \ell_+$ satisfying
$ \ell_-<\lambda_1 < \lambda_3< \ell_+< \lambda_3+\epsilon$
the equation (\ref{eq:PVCnaolinear}) has {\it exactly} $6$ solutions for $g = - t \phi_1$ for  large, positive $t$.
For concreteness, we take $f$ such that $f'(x)=\alpha \ \text{arctan}(x)+\beta$, where $\alpha$ and $\beta$ are adjusted so that $\ell_- = -1$, $\ell_+ = 25.3397$. Finally, set $g =-1000 \ \phi_1 $.

\bigskip
We first obtain a solution $P_0$ by a continuation method.
Set
$ r=\{P_0+s(0.8\phi_1  - 0.1\phi_2  - 0.1\phi_3 )$, a stretch of which ($s \in [-1000, 1000]$) we traverse with an increment $h_s=0.1$. As in the previous section, small terms are inserted so as to increase the possibilities that intersections with the critical set are transversal.

\begin{figure} 
\centering
\begin{minipage}{0.45\textwidth}
	\centering
	\includegraphics[width=165pt]{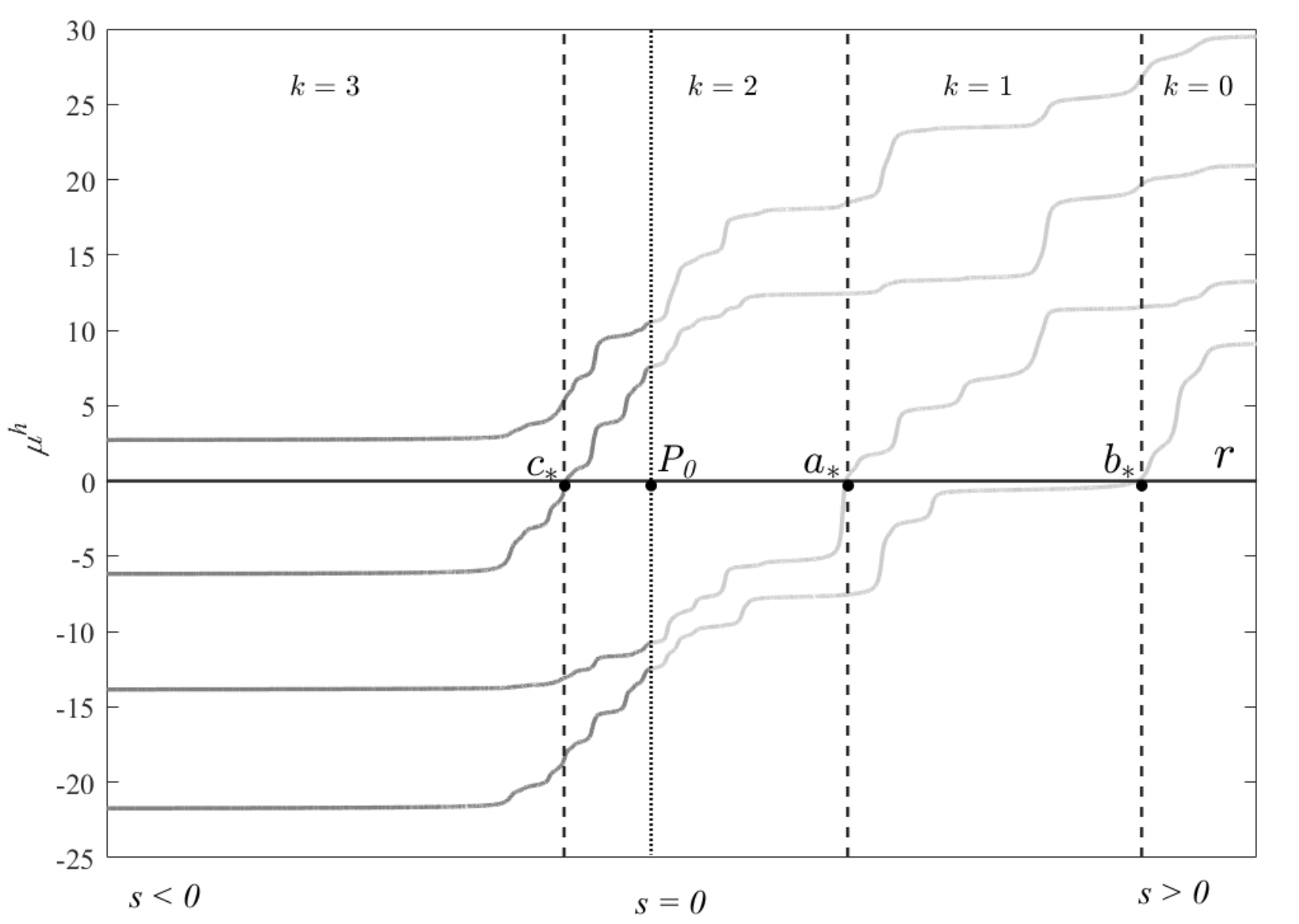}
	\caption{The four smallest eigenvalues of $DF(u)$ for $u(s) \in r$.}
	\label{fig:5.2}
\end{minipage}
\hspace{0.4cm}
\begin{minipage}{0.45\textwidth}
	\centering
	\includegraphics[width=150pt]{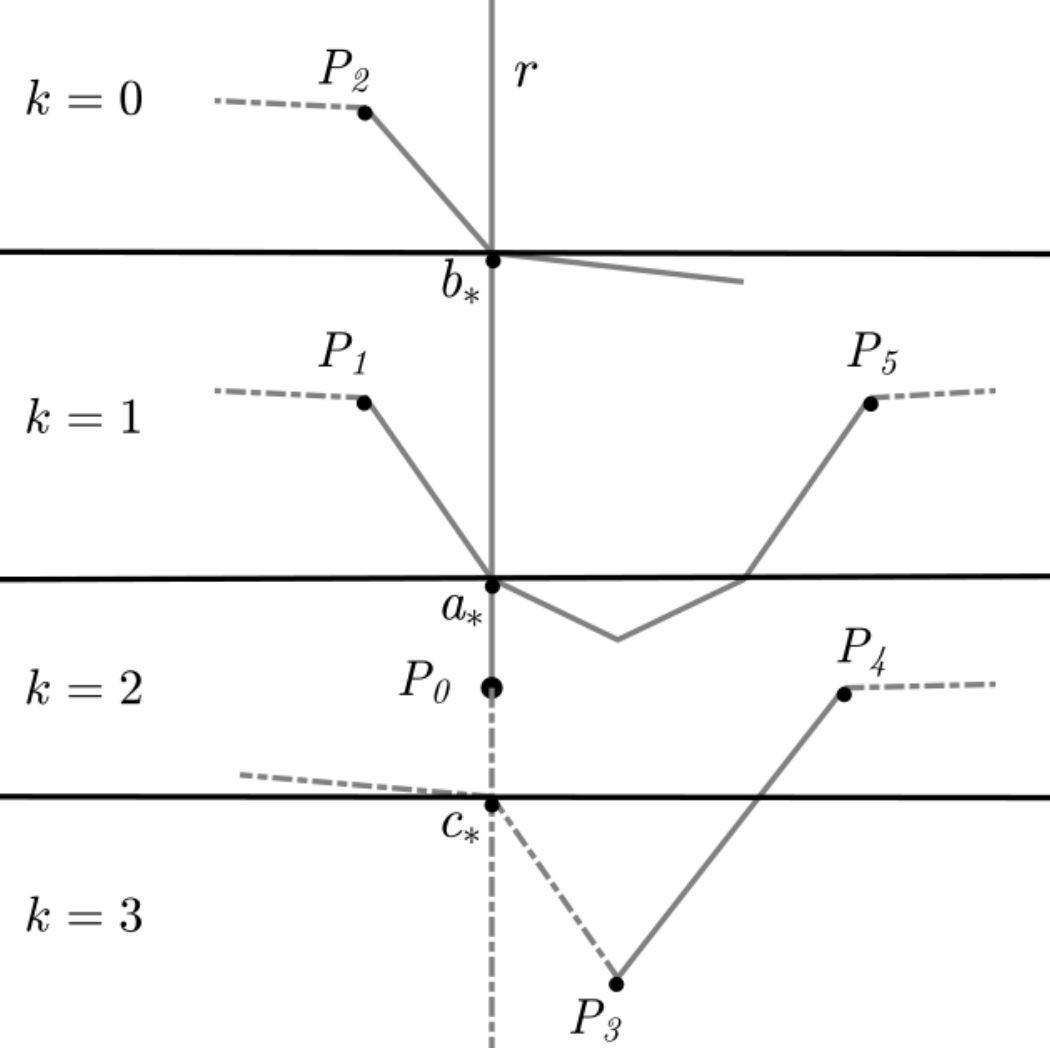}
	\caption{The six solutions in $\mathcal{B}$}
	\label{fig:morse}
\end{minipage}
\end{figure}

Along $r$, the four smallest eigenvalues $\mu_i $ of the Jacobian $DF(u)$ are given in Figure \ref{fig:5.2} (for the underlying numerics, we used \cite{BOFFI}). A point $u \in r$ is a critical point of $F$ if and only if some such eigenvalue is zero.

\begin{figure} 
$$\begin{array}{cc}	
	\includegraphics[keepaspectratio,scale=0.13]{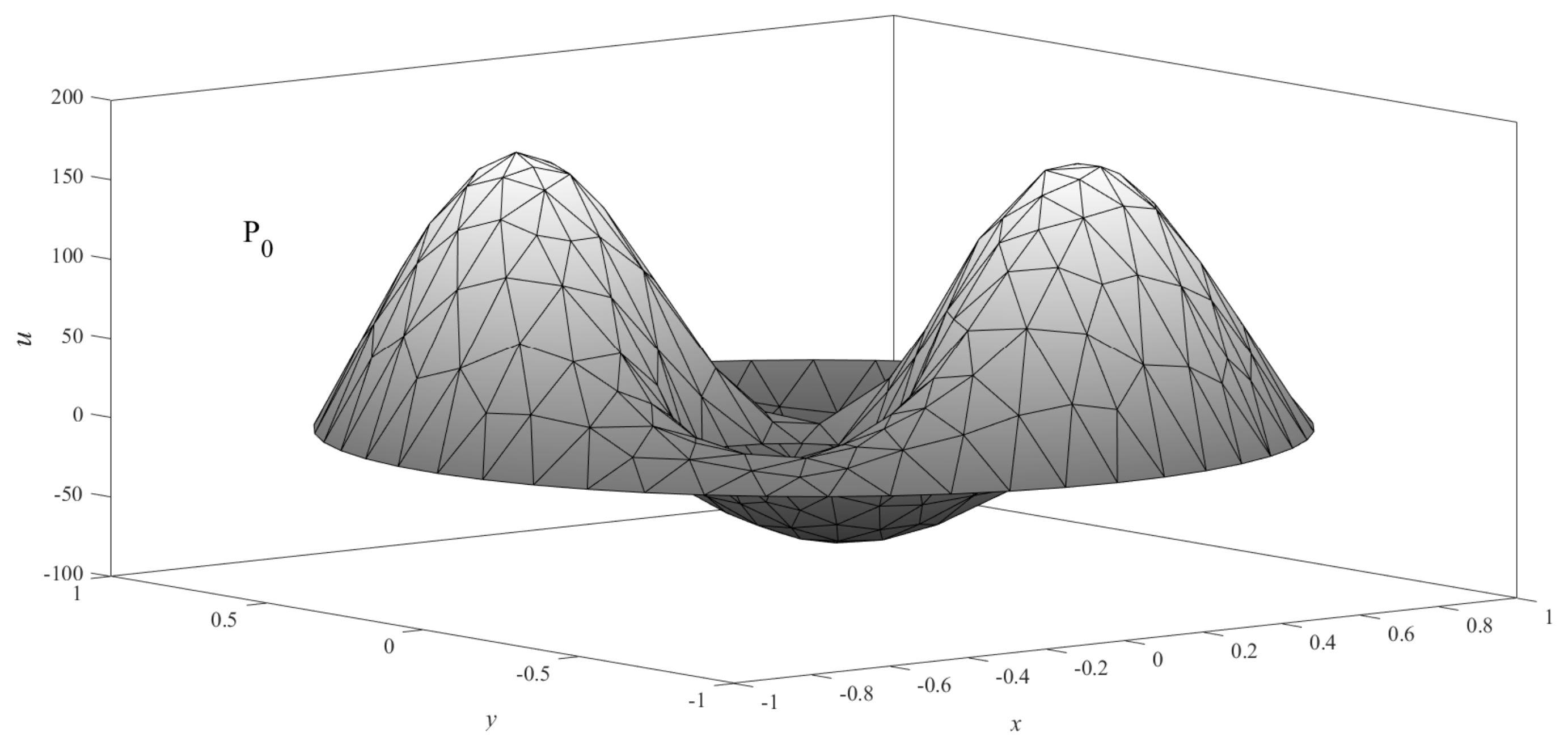} &
	\includegraphics[keepaspectratio,scale=0.13]{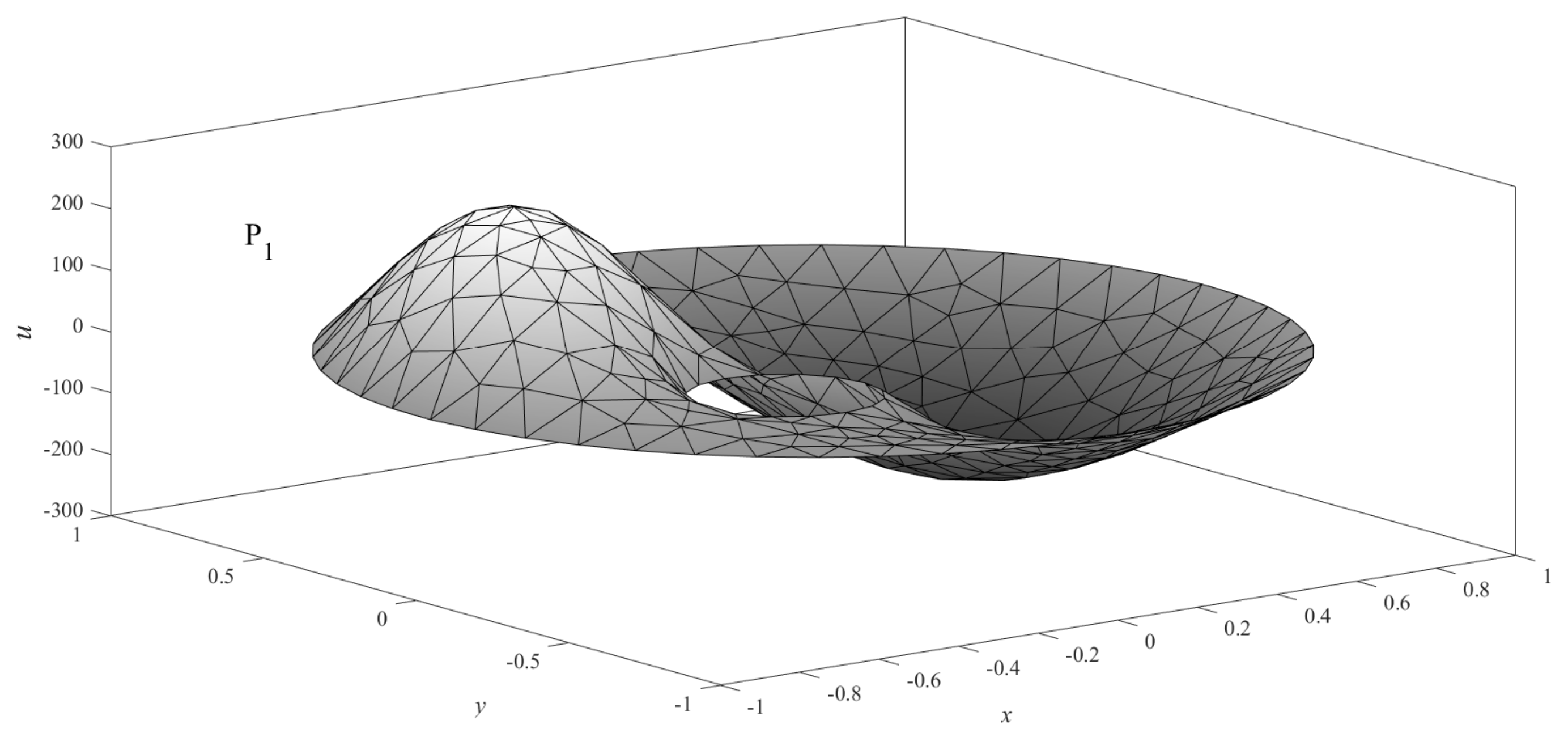} \\
	\includegraphics[keepaspectratio,scale=0.13]{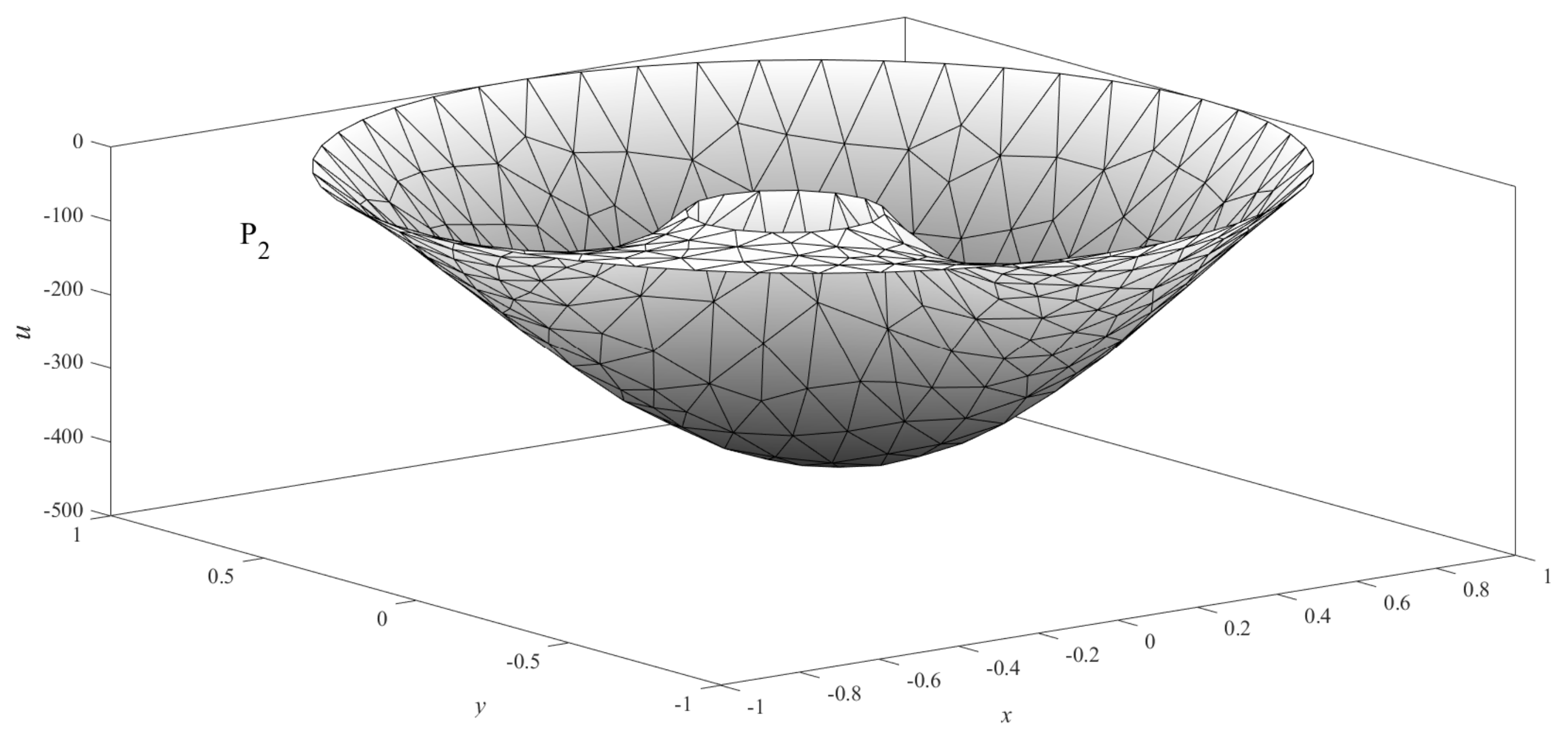} &
	\includegraphics[keepaspectratio,scale=0.13]{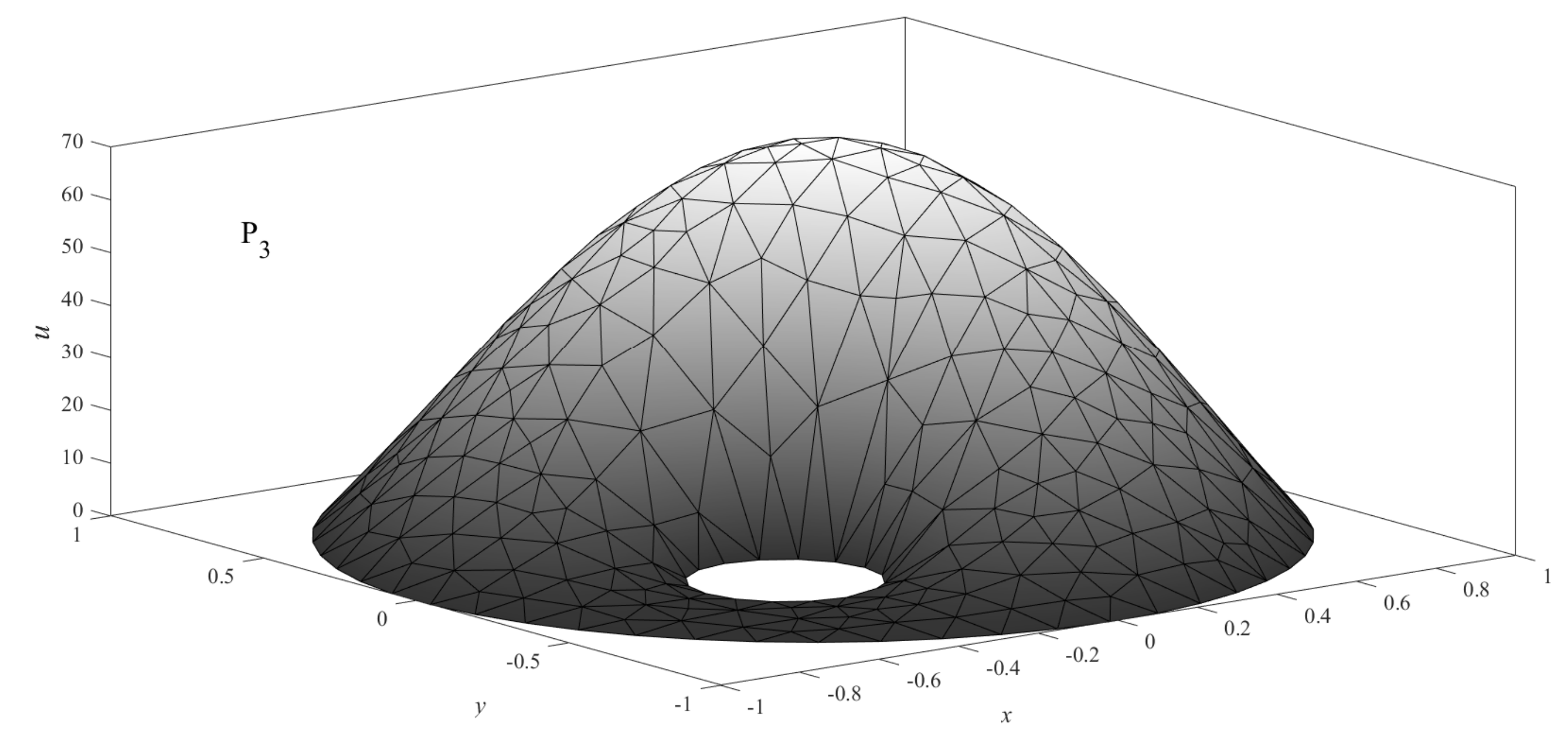} \\
	\includegraphics[keepaspectratio,scale=0.13]{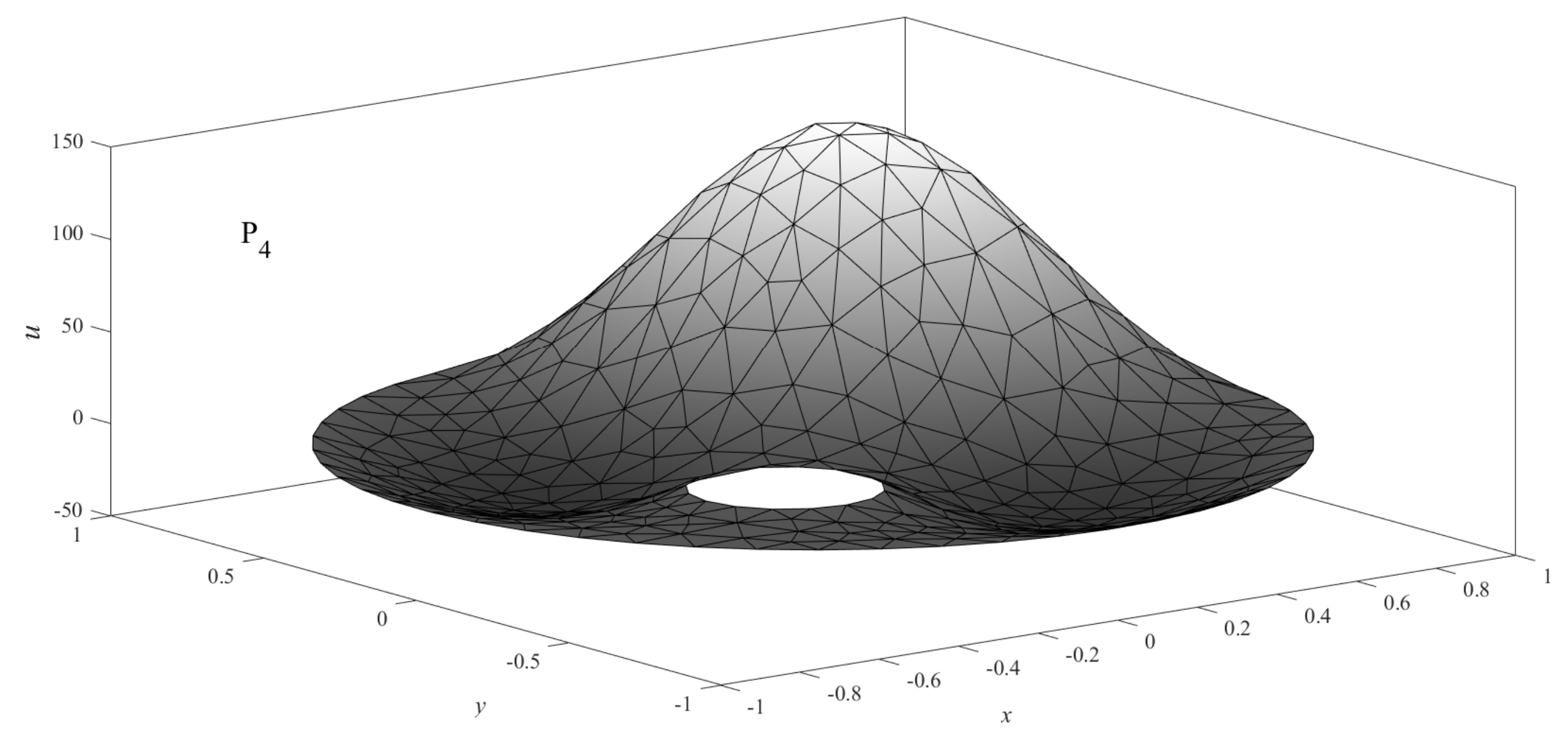} &
	\includegraphics[keepaspectratio,scale=0.13]{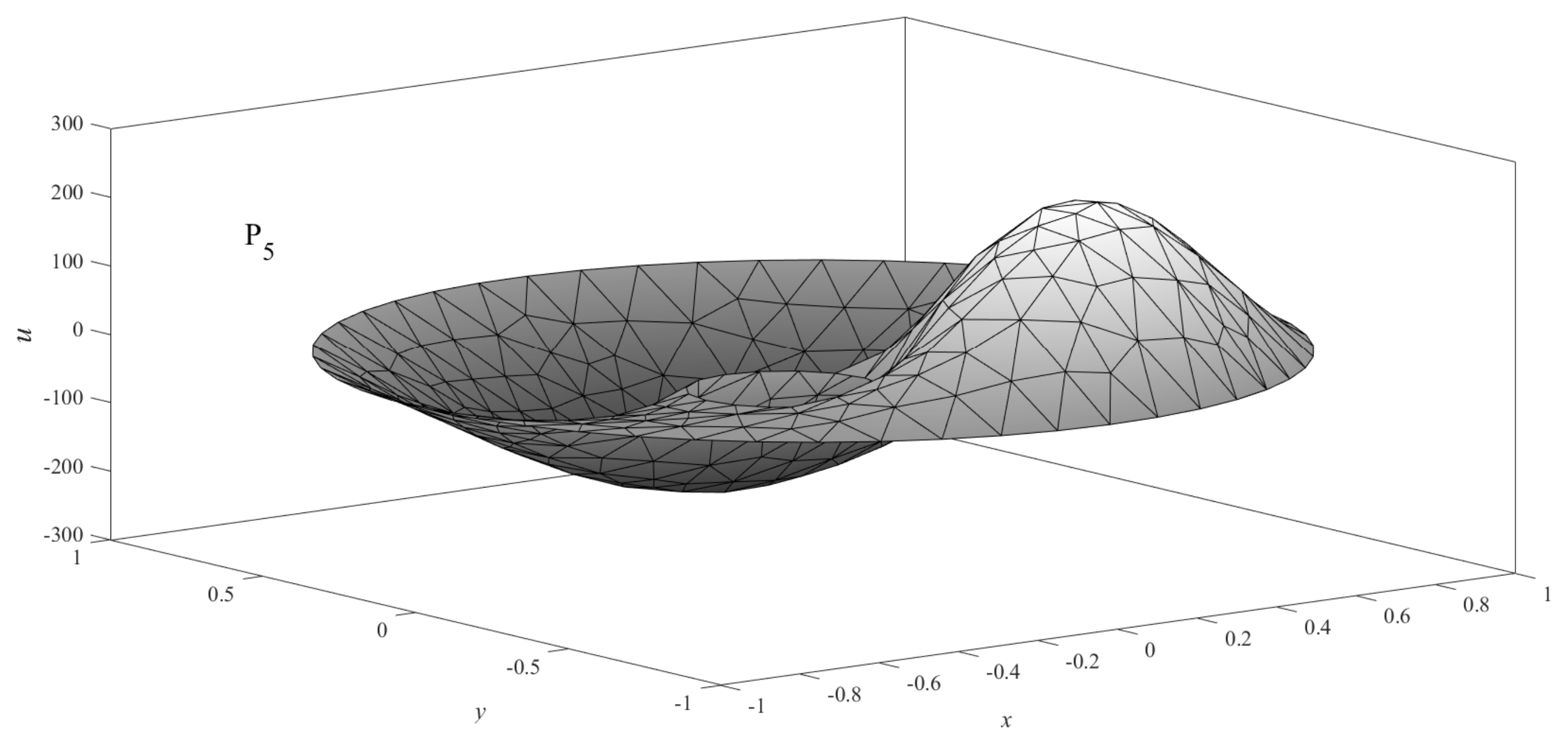}
\end{array}$$
\caption{The six solutions}
\label{fig:5.5}
\end{figure}

In Figure \ref{fig:morse}, horizontal lines represent parts of the critical set $\mathcal{C}$. The value $k$  counts the number of negative eigenvalues of $DF(u)$ at each of the regular components. The point $P_0 \in r$ belongs to a component for which $k=2$. The bifurcation diagram $\mathcal{B}$, containing the six solutions, is described in Figure \ref{fig:morse} and the solutions are given in Figure \ref{fig:5.5}. Continuation to $P_5$ required  finer
jumps along $r$.

No relative residue
$\epsilon({u })=\frac{||F (u )-g ||_{Y }}{||g ||_{Y }}$ is larger than $10^{-12}$.

\appendix
	\section{Some terminology and basic facts} \label{terminology}
	
	We provide definitions and basic properties of some of the concepts.

	\subsection{Fredholm operators} \label{Fredholm}
	
	A linear, bounded operator $L:B_1 \to B_2$ between real, Banach spaces  is a {\it Fredholm operator} if
	$\dim \ker L < \infty$,
	and its range  is closed, admitting a complementary subspace $V$ of finite dimension.
	In a nutshell, being Fredholm  is what it takes for an operator to look like a (rectangular) matrix acting on finite dimensional spaces.
	The {\it index} of a Fredholm operator is $\ind L =  \dim \ker L -\dim V$. Square matrices are operators of index zero.
	As usual, ${\mathcal{B}}(B_1, B_2)$ is the set of bounded operators between the Banach spaces $B_1$ and $B_2$ equipped with the operator norm. We use the following properties  (\cite{LAXfun}).
	
	\medskip
	\begin{enumerate}
		\item Fredholm operators form an open, dense space of ${\mathcal{B}}(B_1, B_2)$.
		\item Composition of Fredholm operators is Fredholm, $\ind L_1 L_2 = \ind L_1 + \ind L_2$.
		\item Let $L$ be Fredholm. If $P \in{\mathcal{B}}(B_1, B_2)$ is sufficiently small or compact, then $L+P$ is also Fredholm and $\ind (L+P) = \ind L$.
	\end{enumerate}
	
	\subsection{Identifying folds} \label{folds}
	
	Let $X$ and $Y$ be real Banach spaces admitting a bounded inclusion $X \hookrightarrow Y$. Denote by $\cB(X,Y)$  the set of bounded operators between $X$ and $Y$, endowed with the usual operator norm. 
	
	\medskip
	
	\begin{theorem} \label{trespassing}
		
		\noindent Suppose that the function $F: U \to Y$ is of class $C^3$. For a fixed $x \in U$, let the Jacobian $DF(x): X \to Y$ be a Fredholm operator of index zero for which 0 is an eigenvalue. Let $\Ker DF(x)$ be spanned by a vector $k \in X$ such that $k \notin \Ran DF(x)$. Then the following facts hold.
		
		\noindent (1) For some open ball $B \subset  \cB(X,Y)$ centered in $DF(x)$, operators $T \in B$ are also Fredholm of index zero and $\dim \Ker T \le 1$. 
		
		\noindent (2) There is a  function $\lambda_s: T \in B \to \RR$ taking $T$ to its eigenvalue of smallest module, which is necessarily a real eigenvalue. The function is real analytic. For a suitable normalization, the  corresponding eigenvector function $T \mapsto \phi_s(T) \in X$ is  real analytic. 
		
		\noindent (3) For a small open ball $B_x \subset X$ centered in $x$, there are $C^3$ functions
		\[ \tilde x \in B_x \mapsto \lambda_s( DF(\tilde x)) \in \RR \ , \quad x \mapsto \phi_s( DF(\tilde x)) \in X\ . \]
		
		If additionally  
		$D \lambda_s (x) . \phi(x) \ \ne \ 0 $,  $x$ is a fold. 
		
	\end{theorem}
	
	\medskip
	In particular, the critical set $\cC$ of $F$ is a submanifold of $X$ of codimension one near $x$. Jacobians are not required to be self-adjoint operators. Characterizations of folds  in the infinite dimensional context may be found in \cite{BCT1, BCT2, D, BD, MST3}.
	
	\begin{proof}  We barely sketch a lengthy argument. The implicit function theorem yield items (1) and (2), as described in Proposition 16 of \cite{CaTZ}. Item (3) then follows. Proposition 2.1 of \cite{MST3} implies item (4).	\end{proof}

	\subsection{Continuation at a fold from spectral data} \label{continuationatafold}
	
	Predictor-corrector methods at regular points are well described in the literature (\cite{ALLGOWER,KELLER,KELLEY,RHEINBOLDT}). Here we provide details about the inversion algorithm we employ in the examples in Section  \ref{Solimini} in the neighborhood of a fold.
	
	The algorithm must identify critical points $u$ of $F:X \to Y$. Due to the nature of the examples, this is accomplished by checking if some eigenvalue of the Jacobian $DF(u)$ is zero. We assume that the original problem $F(u) = g$ admits a variational formulation, so that $DF(u)$ is a self-adjoint operator, and the task is simpler. The general case may also be handled, but we give no details.
	
	We modify the prediction phase of the usual continuation method and perform correction in a standard fashion.
	Following (\cite{UECKER, CALNETO,KAMINSKI}), we use spectral data: by continuity, for $u$ close to a fold $u_c$, the Jacobian $DF(u)$ has an eigenvalue $\lambda$ close to a zero eigenvalue $\lambda_c=0$ of $DF(u_c)$, and a normalized eigenvector $\phi$ close to $\phi_c$, a normalized generator of $\ker DF(u_c)$. The eigenvalue $\lambda$  plays the role of arc length in familiar algorithms.
	
	For a smooth function $F: X \to Y$ between real Banach spaces, we search for the preimage $u(t)$ of a smooth curve $ \gamma(t) \subset Y$  such that, at $t = t_c$, $\gamma(t_c) = F(u_c)$ is the image of a fold $u_c$. As usual, we consider the homotopy
	$$
	\begin{array}{lrll}
		H:&X \times \mathbb{R} &\longrightarrow Y, \ \ (u,t)&\longmapsto F(u)- \gamma(t) \\
	\end{array}
	$$
	and assume the hypothesis of the implicit function theorem:
	$$
	\begin{array}{lrll}
		DH(u,t): X \times  \mathbb{R} \to  Y \ , \quad (\hat u, \hat t) \ \mapsto DF(u) \ \hat u - \gamma' (t) \ \hat t
	\end{array}
	$$
	is surjective at $(u_c, t_c)$. Clearly, $\gamma'(t) \in Y$.
	
	\medskip
	\begin{proposition}  $DH(u_c, t_c)$ is  surjective if and only if $\gamma'(t_c) \notin \Ran DF(u_c)$.
	\end{proposition}
	
	\medskip
	Geometrically, the curve $\gamma(t) \in Y$ crosses the image of the critical set $F(\mathcal{C})$ transversally at the point $\gamma(t_c) = F(u_c)$. Since  $\gamma$ is chosen by the programmer, this is no real restriction.

	\begin{proof}
		As $u_c$ is a fold, $DF(u_c)$ is  a Fredholm operator of index 0 with one dimensional kernel, and  image given by a closed subspace of codimension one. Surjectivity of $DH(u_c, t_c)(\hat u , \hat t) = DF(u_c) \hat u - \gamma'(t_c) \hat t$ holds exactly if $\gamma'(t)$ generates a complementary subspace to $\Ran DF(u_c)$.
	\end{proof}
	
	The next proposition ensures that the inversion of appropriate operators may be performed as in the finite dimensional case.  If $ X = Y = \RR^n$, $DH(z)$ is an $n \times (n+1)$ matrix of rank $n$.
	
	\medskip
	\begin{proposition}
		\label{prop:JacobianaFredholm} For $(u,t)$ close to $(u_c, t_c)$, the Jacobian
		$DH(u,t): X \times \RR \to Y$ is a Fredholm operator of index 1. If $\gamma'(t_c) \ne 0$, $\dim \ker DH(u,t)= 1$ if and only if $u = u_c$, otherwise it is trivial.
	\end{proposition}

	\begin{proof} We first show that $DH(z_c)$ is a Fredholm operator of index 1 with one dimensional kernel.
		Identify $\RR \sim \{  \gamma'(t_c) \hat t , \hat t \in \RR\}$ and write $DH(z_c)$ as the composition
		\[( \hat u , \hat t) \in X \times \RR \mapsto (DF(u_c) \hat u , \gamma'(t_c) \hat t) \in Y \times \RR \ni  (y, s)   \mapsto y - s \in Y \ , \]
		easily seen to consist of Fredholm operators of indices $0$ and $1$ respectively. From  property (2) in Section \ref{Fredholm}, $DH(u_c, t_c)$ is a Fredholm operator of index 1. If $\gamma'(t_c) \ne 0$, then $\dim \ker DH(u_c, t_c) \le 1$ and $\dim \ker DF(u_c) = 1$ if and only if $\dim \ker DH(u_c, t_c)=1$. Property (3) in  Section \ref{Fredholm} also implies that $DH(u,t)$ is Fredholm of index 1, for $(u,t)$ near $(u_c, t_c)$, as $H$ is of class $C^1$. Smoothness of  eigenvalues and eigenvectors proves the claims for $\ker DH(u,t)$.
	\end{proof}
	
	To obtain a prediction from point $(u,t) = (u(t),t) \in H^{-1}(0)$, we must find a nonzero tangent vector $(\hat u, \hat t) \in T_{(u,t)} H^{-1}(0)$, so that
	$$
	\begin{array}{lrll}
		DH(u,t)(\hat u, \hat t)=DF(u) \hat u - \gamma'(t) \hat t =0 \ , \quad (\hat u, \hat t) \ne 0 .
	\end{array}
	$$
	
	At points $u= u(t)$ for which $DF(u)$ is invertible, this is easy: set $\hat t = 1$ and get $\hat u $ by solving a linear system. Instead, we assume $u$ close to a fold $u_c \in \mathcal{C}$. By the smoothness of simple eigenvalues and associated (normalized) eigenvectors,  $DF(u)$ has an eigenvalue $\lambda$ and associated eigenvector $\phi$  near an eigenvalue  $\lambda_c = 0$ and eigenvector  $\phi_c$ of $DF(u_c)$. We must compute a nonzero solution $(\hat u ,\hat t)$ of $DH(u,t)(\hat u ,\hat t) = 0$, or equivalently $DF(u)\hat u = - \gamma'(t) \hat t$,  with a procedure which is continuous in $t \sim t_c$.

	For $\hat u \in X$, split $ \hat u= \hat v + \hat r \phi$ for $\langle \hat v,  \phi \rangle =0$ and $\langle \phi, \phi \rangle =1$. Clearly, $\hat v$ and $\hat r$ are continuous in $u$, since $\phi= \phi(u)$ is. The tensor product $\phi \otimes \phi$ denotes the rank one linear operator $(\phi \otimes \phi) v = \langle \phi, v \rangle \phi$. In particular, as $\| \phi \|=1$, $(\phi \otimes \phi) \phi = \phi$.
	
	\medskip
	\begin{proposition}
		Let $u$ be sufficiently close to a fold $u_c \in X$  of $F$, with the eigenvalue $\lambda \sim 0$ such that $| \lambda|< 1$ and associated normalized eigenvector $\phi$. For $\alpha \ge 1$ the operator $S(u)=DF(u)+\alpha \ \phi \otimes \phi  : X \to Y$ is invertible.	
	\end{proposition}
	
	Notice that $S = DF(u)$ when restricted to $\phi^\perp$.
	\medskip
	
	\begin{proof} The operator $S(u)$ is a rank one perturbation of $DF(u)$, so, by property (3) in  Section \ref{Fredholm}, it is also a Fredholm operator of index zero: invertibility is equivalent to injectivity. For $\hat u= \hat v + \hat r \phi$, with $\hat v \in \phi^\perp$ and $\langle \phi, \phi \rangle =1 $,
		\[ S(u) \hat u= DF(u) (\hat v + \hat r \phi) + \alpha(\phi \otimes \phi) (\hat v + \hat r \phi) = 0 \ , \]
		implies
		\[ S(u) \hat u= DF(u) \ \hat v  +  \hat r \lambda \phi + \alpha \ \hat r \ \phi  = 0 \ . \]
		Since $DF(u)$ is self-adjoint, $\phi^\perp$ is an invariant subspace and both terms are zero, 
		\[ DF(u) \hat v =0 ,\  (\alpha + \lambda) \hat r \ \phi  = 0 \ . \]
		As the restriction $DF(u)$ to $\ker DF(u)^\perp$ is an isomorphism,  $\hat v = 0$, $\hat r = 0 $.
	\end{proof}
	
	\begin{proposition} Under the hypotheses of the proposition above, the solution of
		\[ S(u) \hat u = - \lambda \gamma'(t)  - \alpha \langle \phi, \gamma'(t) \rangle \phi \]
		is of the form $\hat u = \hat v - \langle \phi, \gamma'(t) \rangle \phi$ for some $\hat v \in \phi^\perp$, $\| \phi \| = 1$. Moreover,
		\[ DF(u) \hat u = - \lambda \gamma'(t) \ \hbox{and} \ (\hat u , \lambda) \in T_{(u,t)} H^{-1}(0) \ . \]
	\end{proposition}
	
	In other words, $(\hat u , \lambda)$ is the tangent vector required in the prediction phase.

	\begin{proof} As $S(u)$ is invertible, $\hat u$ is well defined for the given right hand side. For $\hat u = \hat v + \hat r \phi$ with $\hat v \in \phi^\perp$ and $\langle \phi, \phi \rangle =1 $,
		take the inner product with $\phi$ of
		\[ S(u) \hat u = DF(u) \hat u + \alpha \hat r \phi = - \lambda \gamma'(t)  - \alpha \langle \phi, \gamma'(t) \rangle \phi \ . \]
		to obtain $\hat r = - \langle \phi, \gamma'(t) \rangle$ and then
		$DF(u) \hat u = - \lambda \gamma'(t)$ follows.
	\end{proof}
	
	In the application of Section \ref{Solimini} finite element methods applied to the Laplacian yields the usual sparse matrices. The term $\alpha \phi \otimes \phi$ spoils sparseness, and one has to proceed by inversion through standard techniques associated with rank one perturbations. The numerical inversion worked well with $\alpha=1$.


\begin{thebibliography}{99}

	\bibitem{ALLGOWER} {\sc E. L. Allgower and K. Georg}, {\em Numerical continuation  methods: an introduction}, Springer--Verlag, New York, 1991.

\bibitem{ALLGOWER2} {\sc E. L. Allgower, S.  Cruceanu and S. Tavener}, {\em Application of numerical continuation to compute all solutions of semilinear elliptic equations}, Adv.Geom. 9, 2009,  371--400.

\bibitem{AP} {\sc A. Ambrosetti and G. Prodi}, {\em On the inversion of some differentiable mappings with singularities between Banach spaces}, Ann. Mat. Pura Appl. 4, 93, 1972,  231--246.

\bibitem{ARDILA} {\sc L.A. Ardila}, {\em Morin singularities of the McKean-Scovel operator}, Ph.D. Thesis, PUC--Rio, Rio de Janeiro, 2021.

\bibitem{BD}
\textsc{F. Balboni and F. Donati},
Singularities of Fredholm maps with one-dimensional kernels, I: A complete classification,
\textit{arXiv: Functional Analysis} 1,   1--67, 2014.

\bibitem{BERGER} {\sc M. S. Berger and E. Podolak}, {\em On the solutions of a nonlinear Dirichlet
	problem }, Indiana Univ. Math. J., 24, 1974,  837--846.







\bibitem{BCT1}
\textsc{M. S. Berger,  P. T. Church and J. G. Timourian},
Folds and Cusps in Banach Spaces, with Applications to Nonlinear Partial Differential Equations. I,
\textit{Indiana Univ. Math. Journal} 34,   1--19, 1985.

\bibitem{BCT2}
\textsc{M. S. Berger,  P. T. Church and J. G. Timourian},
Folds and Cusps in Banach Spaces, with Applications to Nonlinear Partial Differential Equations. II,
\textit{Transactions of the Amer. Math. Soc.} 307,   225--244, 1988.

\bibitem{BOFFI} {\sc D. Boffi}, {\em Finite element approximation of eigenvalue problems}, Acta Numer. 2010, 1--120.




\bibitem{PLUM} {\sc B. Breuer, P.J.McKenna and M. Plum}, {\em Multiple solutions for a semilinear boundary value problem: a computational multiplicity proof}, J. Diff. Eqs. 195, 2003, 243--269.

\bibitem{BROWN}  {\sc  K. M. Brown and W. B. Gearhart}, {\em Deflation techniques for the calculation of further solutions of a nonlinear system}, Numer. Math. 16 (1971),
334--342.


\bibitem{TOMEIEBUENO} {\sc H. Bueno and C. Tomei}, {\em Critical sets of nonlinear Sturm--Liouville problems of Ambrosetti--Prodi type}, Nonlinearity, 15, 2002, 1073--1077.


\bibitem{BURGHELEA} {\sc D. Burghelea, N.C. Saldanha and C. Tomei}, {\em  Results on infinite--dimensional topology and applications to the structure of the critical set of nonlinear Sturm Liouville operators}, J. Diff. Eqs. 188, 2003, 569--590.

\bibitem{BURGHELEA2} {\sc D. Burghelea, N.C. Saldanha and C. Tomei}, {\em  The topology of the monodromy map of a second order ODE}, J. Diff. Eqs. 227, 2006, 581 - 597.


\bibitem{CALNETO} {\sc J.T. Cal Neto and C. Tomei}, {\em Numerical analysis of semilinear elliptic equations with finite spectral interaction}, J.Math.Anal.Appl., 395, 2012,  63--77.

\bibitem{CaTZ} {\sc M. Calanchi, C. Tomei and A. Zaccur}, {\em Cusps and a converse to the Ambrosetti-Prodi theorem}, Ann. Sc. Norm. Sup. Pisa  XVIII  (2018) p. 483--507.

\bibitem{CHIAPPINELLI} \textsc{R. Chiappinelli} and \textsc{R. Nugari}, {The Nemitskii operator in H\"older spaces: some necessary and sufficient conditions},
J. London Math. Soc. 51 (1995),  365-- 372.




\bibitem{CIARLET} {\sc P.G. Ciarlet}, {\em The finite element method for elliptic problems}, SIAM, Philadelphia, PA, 2002.


\bibitem{COSTA} {\sc D.G. Costa, D.G. Figueiredo and P.N. Srikanth}, {\em The exact number of solutions for a class of ordinary differential equations through Morse index computation},
J. Diff. Eqs. 96, 1992, 185--199.

\bibitem{DUCZMAL} {\sc L. Duczmal}, {\em Geometria e invers\~ao num\'erica de fun\c c\~oes de uma regi\~ao limitada do plano no plano}, Ph.D. Thesis, PUC--Rio, Rio de Janeiro, 1997.

\bibitem{D} {\sc J. Damon}, {\em A Theorem of Mather and the local structure of nonlinear Fredholm maps}, Proc. Sym. Pure Math. 45 part I, (1986) 339-352.

\bibitem{KAMINSKI} {\sc O. Kaminski}, {\em An\'alise Num\'erica de Operadores El\'ipticos Semi--Lineares com Intera\c c\~ao Espectral Finita}, Ph.D. Thesis, PUC--Rio, Rio de Janeiro, 2016.

\bibitem{FARRELL}{\sc P.E. Farrell, \'A. Birkisson, and S.W. Funke}, {Deflation techniques for finding distinct solutions of nonlinear partial differential equations}, SIAM J. Comput. 37(4), A2026--A2045, 2015.


\bibitem{HATCHER} {\sc A. Hatcher}, {\em Algebraic Topology}, Cambridge U. Press, 2002.


\bibitem{KELLER} {\sc H.B. Keller}, {\em Lectures on Numerical Methods in Bifurcation Theory}, Tata Institute of Fundamental Research, Lectures on Mathematics and Physics, Springer, New York, 1987.

\bibitem{KELLEY} {\sc C.T. Kelley}, {\em Numerical methods for nonlinear equations}, Acta Numer. 27, 2018,  207--287.

\bibitem{LAXfun} {\sc P. Lax}, {\em Functional Analysis}, Wiley Pure and Applied Mathematics (2002).

\bibitem{LMCKENNA} {\sc A.C. Lazer and P.J. McKenna}, {\em On a conjecture related to the number of solutions of a nonlinear Dirichlet problem}, Proc. Royal Soc. Edinburgh  95A, 1983, 275--283.





\bibitem{LAZERANDMCKENNA} {\sc A.C. Lazer and P.J. McKenna}, {\em A Symmetry Theorem and Applications to Nonlinear
	Partial Differential Equations}, J. Diff. Eqs. 72, 1988, 95--106.




\bibitem{MST1} {\sc I. Malta, N.C. Saldanha and C. Tomei}, {\em The numerical inversion of functions from the plane to the plane},
Math. Comp. 65, 1996, 1531--1552.



\bibitem{MST3} {\sc I. Malta, N.C. Saldanha and C. Tomei}, {\em Morin singularities and global geometry in a class of ordinary differential operators}, Topol. Meth. Nonlinear Anal., 10, 1997, 137--169.

\bibitem{MM}  {\sc  A. Manes and A.M. Micheletti}, {\em Un'estensione della teoria variazionale classica degli autovalori per operatori ellittici del secondo ordine}, Boll. U. Mat. Ital. 7 (1973) 285--301.

\bibitem{DIEGO} {\sc D. S. Monteiro}, {\em Um m\'etodo de continua\c{c}\~{a}o estruturado para problemas com m\'{u}ltiplas solu\c{c}\~{o}es}, Ph.D. Thesis, PUC--Rio, Rio de Janeiro, 2021.

\bibitem{RS} {\sc M. Reed and B. Simon}, {\em Methods of modern mathematical physics IV, Analysis of Operators}, Academic Press,  1978. 

\bibitem{RHEINBOLDT} {\sc W.C. Rheinboldt}, {\em Numerical continuation methods: a perspective}, J. Comp. A Math. 124, 2000,  229--244.

\bibitem{SOLIMINI} {\sc S. Solimini}, {\em Some remarks on the number of solutions of some nonlinear elliptic problems}, Analyse non lin\'eaire 2, 1985, 143--156.

\bibitem{SMILEY} {\sc M.W. Smiley and C. Chun}, {\em Approximation of the bifurcation function for elliptic boundary value problems}, Numer. Meth. PDE, 16, 2000, 194--213.



\bibitem{TELES} {\sc E. Teles and C. Tomei}, {\em The geometry of finite difference discretizations of semilinear elliptic operators}, Nonlinearity, 25, 2012, 1135--1154.

\bibitem{UECKER}{\sc H. Uecker}, {\em Numerical Continuation and Bifurcation in Nonlinear PDEs}, SIAM, 2021.



\bibitem{Whitney} {\sc H. Whitney}, {\em On singularities of mappings of Euclidean spaces. Mappings of the plane into the plane}, Annals Math. Second Series, 62, 1955, 374--410.








\end{thebibliography}
\end{document}